\documentclass[11pt,reqno]{amsart}

\usepackage[T1]{fontenc}
\usepackage[utf8]{inputenc}

\usepackage{newtxtext}
\usepackage{newtxmath}

\usepackage[
a4paper,
left=27mm,
right=27mm,
top=25mm,
bottom=25mm
]{geometry}

\usepackage{mathtools}   
\usepackage{mathrsfs}
\usepackage{extarrows}

\usepackage[shortlabels]{enumitem}
\usepackage{graphicx}
\usepackage{float}
\usepackage{tikz}
\usetikzlibrary{arrows.meta,patterns,decorations.pathreplacing,calc}
\usepackage{booktabs}
\usepackage{xcolor}
\usepackage{microtype}

\setlist{
	topsep=4pt,
	itemsep=1pt,
	parsep=0pt,
	partopsep=0pt
}

\usepackage{xcolor}

\definecolor{myred}{RGB}{220,0,0}
\definecolor{myblue}{RGB}{30,100,230}
\definecolor{myurlblue}{RGB}{40,110,230}

\usepackage[
colorlinks=true,
linkcolor=myblue,
citecolor=myred,
urlcolor=myurlblue
]{hyperref}
\makeatletter
\renewcommand{\tagform@}[1]{%
	\maketag@@@{\textcolor{myblue}{(\ignorespaces#1\unskip\@@italiccorr)}}%
}
\makeatother
\allowdisplaybreaks[2]

\newtheorem{theorem}{Theorem}[section]
\newtheorem{lemma}[theorem]{Lemma}

\theoremstyle{definition}

\theoremstyle{remark}
\newtheorem{remark}[theorem]{Remark}

\numberwithin{equation}{section}

\newcommand{\R}{\mathbb{R}}

\newcommand{\Z}{\mathbb{Z}}

\author{Yanping Chen$^{*}$}
\thanks{Corresponding author: Yanping Chen. Email: yanpingch@126.com.}
\address{Department of Mathematics, Northeastern University; State Key Laboratory of Synthetical Automation for Process Industries, Shenyang 110004, China}
\email{yanpingch@126.com}

\author{Yongrui Tang}
\address{Department of Mathematics, Northeastern University, Shenyang 110004, China}
\email{yongruit@163.com}

\thanks{There search was partially supported by the National Natural Science Foundation
	of China (Grant No.12525105 and 12371092) and LiaoNing Revitalization Talents Program (Grant No. XLYC2502044)}
\title[]{Sharp variational estimates for Schr\"odinger equation in higher dimensions}
\date{}

\begin{document}

\begin{abstract}

	We study strong $r$-variation estimates for the free
	 Schr\"odinger equation on $\mathbb R^n$, $n\geq3$.
	Combining fractal and weighted restriction estimates with
	mixed space--time estimates, we determine the sharp Sobolev
	regularity, up to endpoints, in several ranges of the spatial
	exponent $p$ and the variation exponent $r$.
	In particular, our results give a sharp $L^2(\mathbb R^n)$
	variational extension  of the
	Schr\"odinger maximal estimate of Du and Zhang
	[\emph{Ann. of Math.}, \textbf{189}, 837--861 (2019)].

\end{abstract}

\maketitle
\tableofcontents
\pagenumbering{arabic}

\section{Introduction}
For initial data $f$, the free  Schr\"odinger equation is given by
$$
u(x,t)=e^{it\Delta}f(x)
=\frac{1}{(2\pi)^n}\int_{\R^n}
e^{i(x\cdot\xi-t|\xi|^2)}\mathcal{F}_x(f)(\xi)\,d\xi.
$$
It satisfies
$$
\left\{\begin{array}{ll}
i\partial_tu(x,t)+\Delta u(x,t)=0,
 &(x,t)\in\R^n\times\R,\\
u(x,0)=f(x),&x\in\R^n.
\end{array}\right.
$$
For a finite interval $I=[a,b]$ and $1\leq r<\infty$, define
$$
V_I^r\{F(x,t)\}_t
=\sup_{a\leq t_0<\cdots<t_J\leq b}
\left(\sum_{j=1}^J|F(x,t_j)-F(x,t_{j-1})|^r\right)^{1/r}.
$$
For $r=\infty$, we set $V_I^\infty\{F(x,t)\}_t=2\sup_{t\in I}|F(x,t)|$.
We write $V_I^r\{F\}$ when the family is understood. For any $t_0\in I$,
$$
\sup_{t\in I}|F(x,t)|\leq |F(x,t_0)|+V_I^r\{F(x,t)\}_t.
$$

In the present paper, we consider
\begin{equation}\label{eq1}
	\|V^r_{[0,1]}\{e^{it\Delta}f\}\|_{L^p(\R^n)}
	\lesssim \|f\|_{H^s(\R^n)},\qquad n\geq3.
\end{equation}
We first give the
notation for the ranges in which we determine the least Sobolev
regularity. Let
$$
A_n(p)=\frac n2-\frac np,\qquad
q_n=\frac{2(n+1)}{n-1},\qquad b_n=\frac{n}{n+1}.
$$
For $3\leq m\leq n+1$, define
$$
 p_{w,m}=\left(\frac12-\frac1{m(m-1)(m-2)}\right)^{-1},\qquad
 \gamma_{w,m}=\frac{m-1}{2m}
 -\frac{n+m-1}{2m(m-1)(m-2)}.
$$
Set $\gamma_3=(2-n)/12$ and, for $4\leq m\leq n+1$, define
$$
 \alpha_m=\frac{m^2-3m}{m^2-2m-2},\qquad
 \gamma_m=\alpha_m\gamma_{m-1}+(1-\alpha_m)\gamma_{w,m}.
$$
Let
$$
 M=\max\{m\in\{3,\ldots,n\}:\gamma_j<0\text{ for }3\leq j\leq m\},
 \qquad a=\frac{\gamma_{w,M+1}}{\gamma_{w,M+1}-\gamma_M},
$$
and define
\begin{equation}\label{eq2}
 d_n=\left(\frac{a(M-1)}{2M}+\frac{1-a}{p_{w,M+1}}\right)^{-1}.
\end{equation}
This is the threshold in the weighted restriction estimate of
Du--Li \cite{DL2025}. For example,
$d_3=96/35$, $d_4=636/245$, and $d_5=3648/1459$.
Let
\begin{equation}\label{eq3}
G_n=\frac{2(n+1)((n+2)d_n-2(n+3))}
{n(n+1)d_n-2(n^2+2n-1)},\qquad
\frac{1}{Q_n}=\frac n2-\frac{n+1}{G_n}.
\end{equation}
For $2\leq p\leq G_n$, let
\begin{equation}\label{eq4}
\theta_p=\frac{1/2-1/p}{1/2-1/G_n},\qquad
\frac1{\rho_n(p)}=\frac{1-\theta_p}{q_n}+\frac{\theta_p}{Q_n}.
\end{equation}
For $d_n\leq p\leq G_n$, let
\begin{equation}\label{eq5}
\frac1{\eta_n(p)}=\frac{1/d_n-1/p}{1/d_n-1/G_n}\frac1{Q_n},
\end{equation}
where $\eta_n(d_n)=\infty$. We show below that
$2<d_n<G_n<q_n<Q_n$.

\begin{theorem}\label{th1}
Let $n\geq3$ and $2\leq p\leq\infty$.  If $s>s(p,r)$, then
\begin{equation*} 
	\|V^r_{[0,1]}\{e^{it\Delta}f\}\|_{L^p(\R^n)}
	\lesssim \|f\|_{H^s(\R^n)}
\end{equation*} 
holds  and fails if $s<s(p,r)$, where $s(p,r)$ is defined as: 
\begin{enumerate}[(i)]
\item If $1\leq r\leq2$, then
$$
s(p,r)=\max\left\{\frac2r+\frac1p-\frac12,\,
A_n(p)+\frac2r-1\right\}.
$$
\item If $2\leq r\leq\infty$ and $p\geq G_n$, then
$$
s(p,r)=\max\left\{A_n(p),\,\frac1p+\frac1r\right\}.
$$
\item If $2\leq p\leq G_n$ and $2\leq r\leq\rho_n(p)$, then
$s(p,r)=1/p+1/r$.
\item If $d_n\leq p\leq G_n$ and $\eta_n(p)\leq r\leq\infty$,
then $s(p,r)=A_n(p)$.
\item If $p=2$ and $q_n\leq r\leq\infty$, then $s(2,r)=b_n$.
\end{enumerate}
\end{theorem}

For the remaining exponents, we have the following upper bounds.

\begin{theorem}\label{th2}
Let $n\geq3$, $2\leq p\leq G_n$, and
$\rho_n(p)\leq r\leq\infty$. Let
\begin{equation}\label{eq6}
 \sigma_n(p,r)=A_n(p)+b_n\max\left\{
 0,\,
 \frac{1/p-1/d_n}{1/2-1/d_n},\,
 \frac{1/p-1/d_n+(Q_n/r)(1/d_n-1/G_n)}
 {1/2-1/d_n+(Q_n/q_n)(1/d_n-1/G_n)}
 \right\}.
\end{equation}
If $s>\sigma_n(p,r)$, then
$$
 \left\|V^r_{[0,1]}\{e^{it\Delta}f\}\right\|_{L^p(\R^n)}
 \lesssim_{p,r,s}\|f\|_{H^s(\R^n)}.
$$

\end{theorem}

The exponent agrees with Theorem \ref{th1} on the boundary:
$$
 \sigma_n(p,\rho_n(p))=\frac1p+\frac1{\rho_n(p)},
 \qquad 2\leq p\leq G_n,
$$
and $\sigma_n(p,r)=A_n(p)$ if
$d_n\leq p\leq G_n$ and $r\geq\eta_n(p)$.

\begin{remark}
	By the definition of $V^\infty$, taking $r=\infty$ in
	Theorem \ref{th1} gives
	\begin{align*}
		\left\|\sup_{0\leq t\leq1}|e^{it\Delta}f|\right\|_{L^p(\R^n)}
		&=\frac12
		\left\|V^\infty_{[0,1]}\{e^{it\Delta}f\}\right\|_{L^p(\R^n)}
		\lesssim \|f\|_{H^s(\R^n)},
	\end{align*}
	provided that
	\[
	\begin{cases}
		s>\displaystyle\frac{n}{n+1},&p=2,\\[6pt]
		s>\displaystyle\frac n2-\frac np,&d_n\leq p\leq\infty.
	\end{cases}
	\]
	In these spatial ranges, the estimate fails below the corresponding
	Sobolev thresholds.
	
	For $p=2$, Rogers' local--global equivalence \cite{R2008}
	shows that the global estimate for every $s>n/(n+1)$ is
	equivalent to the local estimate
	\[
	\left\|\sup_{0<t\leq1}|e^{it\Delta}f|\right\|_{L^2(B^n(0,1))}
	\lesssim_s \|f\|_{H^s(\R^n)},
	\qquad s>\frac{n}{2(n+1)}.
	\]
	This recovers the sharp local $L^2$ maximal estimate of
	Du and Zhang \cite[Theorem 1.2]{DZ2019}, up to the endpoint.
	Thus the index $b_n=n/(n+1)$ is the global counterpart of
	the local index $n/(2(n+1))$.

\end{remark}
\begin{remark}
	The ranges not covered by the sharp estimates in Theorem \ref{th1} are
	$$
	\left\{
	\begin{array}{ll}
		2<p<d_n,&\rho_n(p)<r\leq\infty,\\[4pt]
		d_n\leq p<G_n,&\rho_n(p)<r<\eta_n(p).
	\end{array}
	\right.
	$$
	Theorem \ref{th2} gives upper bounds throughout these ranges.
	For $2<p<G_n$ and $r>\rho_n(p)$,
	$$
	\sigma_n(p,r)<\frac1p+\frac1{\rho_n(p)},
	$$
	which improves the bound obtained directly from \eqref{eq31}
	and Theorem \ref{th4}. Sharpness of the resulting upper bounds
	in the displayed ranges remains open.
\end{remark}

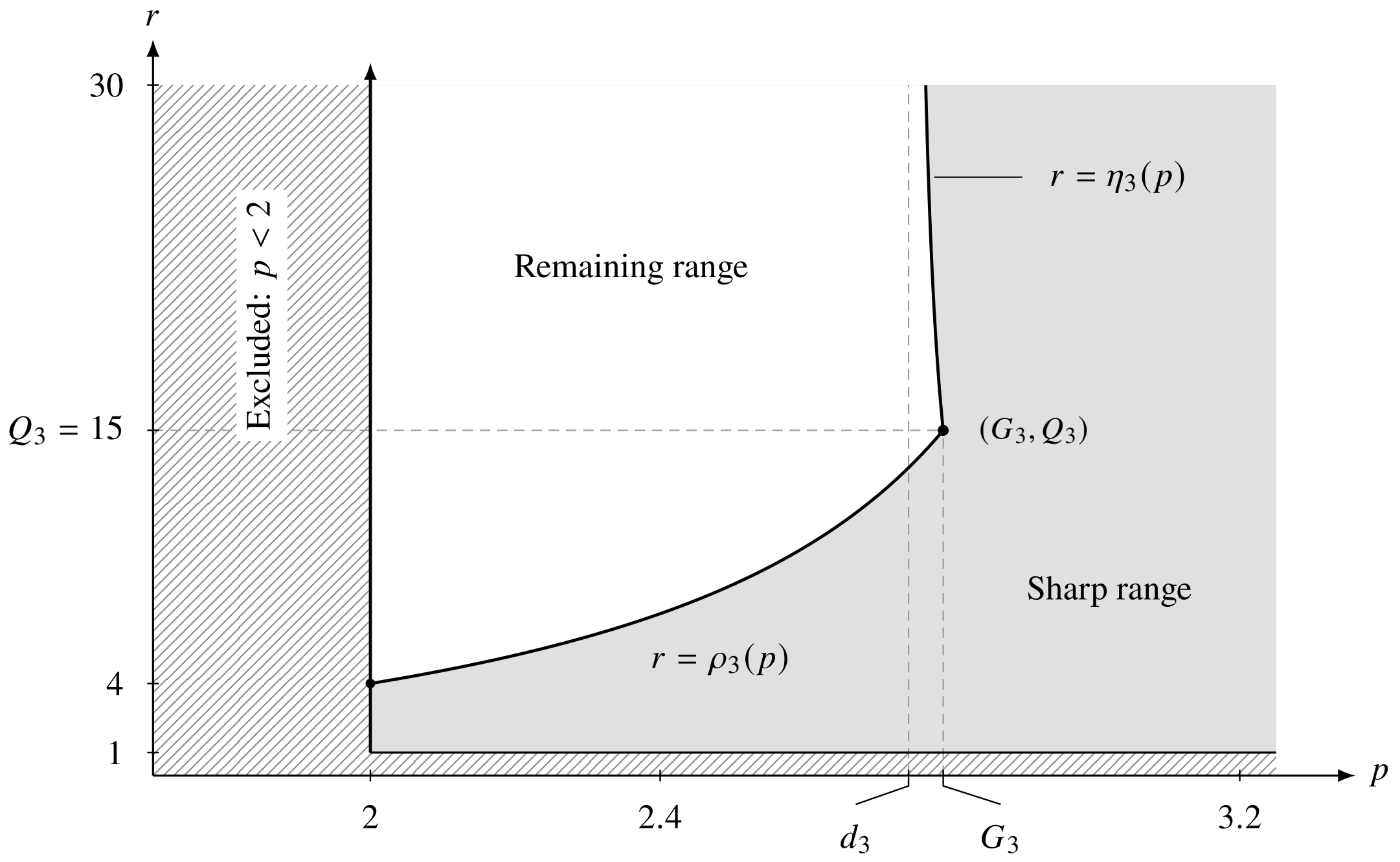
\begin{figure}[H]
	\centering
	\begin{tikzpicture}[x=7.4cm,y=.235cm,
		font=\small,
		axis/.style={black,line width=.6pt,-{Latex[length=1.8mm]}},
		guide/.style={black!40,densely dashed,line width=.4pt},
		boundary/.style={black,line width=.9pt},
		tick/.style={black,line width=.45pt}]
		
		\def\pD{2.742857143}
		\def\pG{2.790697674}
		
		\fill[black!12] (2,1) rectangle (3.25,30);
		
		\fill[white]
		(2,4)
		-- plot[domain=2:\pG,samples=120]
		(\x,{68*\x/(88-27*\x)})
		-- plot[domain=15:30,samples=80]
		({96*\x/(35*\x-9)},\x)
		-- (2,30) -- cycle;
		
		\fill[pattern=north east lines,pattern color=black!45]
		(1.7,0) rectangle (2,30);
		\fill[pattern=north east lines,pattern color=black!45]
		(2,0) rectangle (3.25,1);
		
		\draw[guide] (1.7,15) -- (\pG,15);
		\draw[guide] (\pG,0) -- (\pG,15);
		\draw[guide] (\pD,0) -- (\pD,30);
		
		\draw[boundary]
		plot[domain=2:\pG,samples=120]
		(\x,{68*\x/(88-27*\x)});
		\draw[boundary]
		plot[domain=15:30,samples=80]
		({96*\x/(35*\x-9)},\x);
		
		\draw[boundary,-{Latex[length=1.8mm]}]
		(2,1) -- (2,31);
		\draw[black,line width=.6pt]
		(2,1) -- (3.25,1);
		
		\fill[black] (\pG,15) circle[radius=1.6pt];
		\fill[black] (2,4) circle[radius=1.4pt];
		
		\draw[axis]
		(1.7,0) -- (3.36,0) node[right] {$p$};
		\draw[axis]
		(1.7,0) -- (1.7,32) node[above] {$r$};
		
		\foreach \p/\labeltext in {2/2,2.4/{2.4},3.2/{3.2}}
		{
			\draw[tick] (\p,.25) -- (\p,-.25)
			node[below=3pt] {$\labeltext$};
		}
		
		\foreach \r/\labeltext in {1/1,4/4,15/{Q_3=15},30/30}
		{
			\draw[tick] (1.692,\r) -- (1.708,\r)
			node[pos=0,left=3pt] {$\labeltext$};
		}
		
		\draw[tick]
		(\pD,.25) -- (\pD,-.45)
		-- (2.67,-1.25) node[below=2pt] {$d_3$};
		
		\draw[tick]
		(\pG,.25) -- (\pG,-.45)
		-- (2.87,-1.25) node[below=2pt] {$G_3$};
		
		\node[rotate=90,fill=white,inner sep=3pt]
		at (1.85,20) {Excluded: $p<2$};
		
		\node at (2.36,22) {Remaining range};
		\node at (3.02,8) {Sharp range};
		
		\node[below right=2pt]
		at (2.36,{68*2.36/(88-27*2.36)})
		{$r=\rho_3(p)$};
		
		\node[anchor=west]
		at (2.92,26) {$r=\eta_3(p)$};
		
		\draw[black,line width=.4pt]
		(2.90,26) -- ({96*26/(35*26-9)+.008},26);
		
		\node[anchor=west,font=\footnotesize]
		at (2.82,15) {$(G_3,Q_3)$};
		
	\end{tikzpicture}
	\caption{The ranges in Theorem \ref{th1} for $n=3$,
		shown in the $(p,r)$-plane.
		The white region is the remaining range, with upper bounds
		given by Theorem \ref{th2}.}
	\label{fig1}
\end{figure}

In dimension three, $d_3=96/35$, $G_3=120/43$, $q_3=4$, and $Q_3=15$.
Figure \ref{fig1} shows the ranges in Theorem \ref{th1}
when $n=3$.

The regularity problem for Schr\"odinger variation is closely related to
Carleson's question on pointwise convergence \cite{C1980}: for which $s$
does
\[
\lim_{t\to0}e^{it\Delta}f(x)=f(x)
\quad\text{for almost every }x\in\R^n
\]
hold for every $f\in H^s(\R^n)$? In one dimension, Carleson proved
convergence for $s\geq1/4$, and the examples of Dahlberg and Kenig
\cite{DK1982} exclude $s<1/4$. For $n\geq2$, Sj\"olin \cite{S1987}
and Vega \cite{V1988} independently established convergence for $s>1/2$.
Subsequent improvements connected the problem with Fourier restriction;
see Bourgain \cite{B1995}, Moyua, Vargas, and Vega \cite{MVV1996},
and Tao and Vargas \cite{TV2000}. Lee \cite{L2006} obtained $s>3/8$
in two dimensions using bilinear restriction. Bourgain \cite{B2013}
proved convergence for $s>1/2-1/(4n)$ by multilinear methods and
also showed that the threshold $1/4$ cannot hold in all dimensions.

The necessary regularity was further studied by Luc\`a and Rogers
\cite{LR2017} and by Demeter and Guo \cite{DG2016}, whose argument uses
the pseudoconformal transformation. Bourgain's counterexample
\cite{B2016} gives the condition $s\geq n/(2(n+1))$.
Du, Guth, and Li \cite{DGL2017} reached the corresponding sufficient
condition $s>1/3$ in dimension two. Their proof combines polynomial
partitioning with the decoupling theorem of Bourgain and Demeter
\cite{BD2015}. In higher dimensions, Du and Zhang \cite{DZ2019}
established the local $L^2$ maximal estimate for
$s>n/(2(n+1))$ through a fractal restriction theorem.
Rogers' local--global equivalence \cite{R2008} relates this local
threshold to the global threshold $b_n=n/(n+1)$, as discussed above.

For spatial exponents $p>2$, the maximal problem has a different
dependence on dimension. Rogers and Villarroya \cite{RV2007}
studied estimates over the whole spatial domain. Du, Kim, Wang,
and Zhang \cite{DKWZ2020} obtained further necessary conditions
by using lower-dimensional versions of Bourgain's example.
Cao, Miao, and Wang \cite{CMW2021} and Du and Li \cite{DL2025}
extended the available ranges of higher-dimensional $L^p$ maximal
estimates. The weighted restriction estimate of Du and Li supplies
the exponent $d_n$ used in the present paper.

Variation inequalities also describe repeated oscillations of an
operator family. In particular, finite $r$-variation on $(0,1]$,
for $1\leq r<\infty$, guarantees the existence of the limit at $0$.
The martingale inequalities of L\'epingle \cite{L1976} and
Pisier and Xu \cite{PX1988}, and Bourgain's work on arithmetic
ergodic averages \cite{B1989}, are fundamental examples.
In harmonic analysis, Campbell, Jones, Reinhold, and Wierdl
\cite{CJRW2000} proved oscillation and variation bounds for the
Hilbert transform, while Jones, Seeger, and Wright \cite{JSW2008}
developed variational and jump estimates for several operator families.
Oberlin, Seeger, Tao, Thiele, and Wright \cite{OSTTW2012} obtained
a variation-norm Carleson theorem. For oscillatory singular integrals,
Guo, Roos, and Yung \cite{GRY2020} proved sharp variation estimates
in the modulation parameter for the family studied by Stein and
Wainger \cite{SW2001}.

For the  Schr\"odinger equation, Chen and Tang \cite{CT2026}
determined the sharp Sobolev thresholds for time variation in
two dimensions. In the higher-dimensional problem, the local
estimates must keep track of how often a solution is large.
If $K$ local contributions have size comparable to $\lambda$,
their $\ell^r$ norm has size comparable to $\lambda K^{1/r}$.
The dependence on $K$ is therefore essential when passing from
a maximal estimate to a variation estimate.

Mixed space--time estimates provide a complementary way to control
the time variable. Rogers \cite{R2008} proved local smoothing
estimates and a local--global equivalence, and Rogers and Seeger
\cite{RS2010} obtained endpoint maximal and smoothing estimates.
The work of Lee, Rogers, and Vargas \cite{LRV2011} and Lee,
Rogers, and Seeger \cite{LRS2013} connects mixed-norm bounds
with bilinear and adjoint restriction estimates. Li, Miao,
and Yu \cite{LMY2024} developed this approach on the critical
line $(n+1)/P+1/Q=n/2$. We combine their argument with the
weighted estimate of Du and Li to obtain the mixed-norm input
for the higher-dimensional variation problem.

The proof of Theorems \ref{th1} and \ref{th2} uses the following steps.
Throughout this description, the local estimates concern data with
$\operatorname{supp}\widehat f\subset\{1/2\leq|\xi|\leq2\}$.
\begin{enumerate}[(1)]
\item
We first establish the local $L^2$ variation estimate
\[
\left\|V^{q_n}_{[0,R]}\{e^{it\Delta}f\}\right\|_{L^2(B^n(0,R))}
\lesssim_\varepsilon R^{n/(2(n+1))+\varepsilon}\|f\|_2.
\]
As in the localization argument of \cite{CT2026}, unit space--time
cubes reduce the variation norm to a discrete mixed norm of the
local $L^2$ quantities $f_{j,m}$ defined in \eqref{eq13}.
The cubes are grouped according to their amplitude and the number
$K$ selected above each spatial cube. Summing their squared
contributions gives a factor $K^{1/2}$, whereas the density term
in the Du--Zhang estimate contributes $K^{1/(n+1)}$.
The difference
\[
\frac12-\frac1{n+1}=\frac1{q_n}
\]
determines the variation exponent. A decreasing rearrangement and
a distribution-function argument give the strong discrete
$\ell^{q_n}$ bound, with logarithmic factors absorbed into
$R^\varepsilon$. This yields the estimates at $p=2$ for
$q_n\leq r\leq\infty$.

\item
The second input is a critical mixed-norm estimate with an improved
spatial range. A Galilean transformation converts the weighted
restriction estimate of Du and Li \cite{DL2025} into the shifted-frequency
maximal estimate of Lemma \ref{le8}. We then follow
Li--Miao--Yu \cite[Section 5]{LMY2024}, combining this bound
with Tao's bilinear restriction estimate \cite{Tao2003},
$\varepsilon$-removal, and the bilinear-to-linear argument.
Lemma \ref{th3} gives
\[
\|e^{it\Delta}f\|_{L_x^P L_t^{Q(P)}(\R^{n+1})}
\lesssim_P\|f\|_2,
\qquad
G_n<P\leq q_n,
\qquad
\frac{n+1}{P}+\frac1{Q(P)}=\frac n2.
\]
Thus the weighted maximal estimate enters the variation argument
through the spatial threshold $G_n$ in the critical mixed norm.

\item
For $1\leq r\leq2$, finite-sequence inequalities reduce the local
variation bound to a square-function estimate. Interpolation between
Plancherel's identity and Vega's mixed-norm estimate \cite{V1988},
together with Bernstein's inequality, gives the required spatial
bounds. For $r>2$, we combine the critical mixed-norm estimate,
the $L^2$ variation estimate, and the local maximal estimate.
Interpolation in both $1/p$ and $1/r$ produces the curves
$r=\rho_n(p)$ and $r=\eta_n(p)$ appearing in Theorem \ref{th1}.

\item
Parabolic scaling, spatial localization, and Littlewood--Paley
summation pass from the local estimates to global Sobolev bounds.
The localization follows the methods in
\cite{R2008,LRV2011,LRS2013}, with the variation norm retained
throughout. Theorem \ref{th4} shows that a local factor
$R^{\beta+\varepsilon}$ gives the global condition
\[
s>A_n(p)+2\beta.
\]
For Theorem \ref{th2}, we interpolate the local bounds at
$(2,r_0)$, $q_n\leq r_0\leq\infty$, together with those at
$(G_n,Q_n)$ and $(d_n,\infty)$. The resulting threshold is
$\sigma_n(p,r)$. When $2<p<G_n$ and $r>\rho_n(p)$, it is
strictly smaller than the threshold obtained by keeping the
spatial exponent fixed and extending the bound at $r=\rho_n(p)$.

\item
The necessary conditions are proved by frequency-localized examples.
We adapt the space--time constructions in \cite{LRS2013} and
the variation argument in \cite{GRY2020} to obtain the lower bounds
\[
A_n(p),\qquad
\frac1p+\frac1r,\qquad
\frac2r+\frac1p-\frac12,\qquad
A_n(p)+\frac2r-1.
\]
For $p=2$, rescaling Bourgain's maximal-function example
\cite{B2016,Pierce2020} gives the additional condition $s\geq b_n$.
These examples match the upper bounds in the ranges of
Theorem \ref{th1}. The sharpness of the upper bounds in
Theorem \ref{th2} outside those ranges remains open.
\end{enumerate}

Section \ref{sec1} collects the restriction estimates and interpolation
tools, and derives the shifted-frequency maximal and critical mixed-norm
bounds. Section \ref{sec2} establishes the local variation estimates.
Section \ref{sec4} proves the local-to-global principle and the
sufficiency parts of Theorems \ref{th1} and \ref{th2}.
The counterexamples are given in Section \ref{sec5}.

\section{Preliminaries}\label{sec1}

TThe implicit constants may depend
on $n$ and on the fixed exponents. Dependence on an additional
parameter will be indicated when necessary. We write
$\widehat f=\mathcal F_x(f)$.
The notation $A\lesssim_\varepsilon B$ means
$A\leq C_\varepsilon B$. All estimates are first proved for
Schwartz functions.
In Section \ref{sec2}, $R$ denotes the radius and time
length of the unit-frequency local problem. In
Section \ref{sec4}, the data $g_R$ have frequency comparable
to $R$; parabolic scaling then uses the local estimate at radius $R^2$.

\begin{lemma}\label{le4}
	Suppose $0\leq\theta\leq1$ and
	$$
	\frac1p=\frac{1-\theta}{p_0}+\frac{\theta}{p_1},
	\qquad
	\frac1r=\frac{1-\theta}{r_0}+\frac{\theta}{r_1},
	\qquad 1\leq p_i,r_i\leq\infty.
	$$
	For $r<\infty$,
	\begin{equation}\label{eq12}
		\|V^r_I F\|_p
		\leq \|V^{r_0}_I F\|_{p_0}^{1-\theta}
		\|V^{r_1}_I F\|_{p_1}^{\theta},
	\end{equation}
	with the convention for $V^\infty$ made in the introduction.
	For $1\leq r_0\leq r<\infty$, one also has
	$V^r_I F\leq V^{r_0}_I F$.
\end{lemma}
\begin{proof}
	For a fixed partition, apply the interpolation inequality for
	finite sequences to its increments. Here $\max_j|F(x,t_j)-F(x,t_{j-1})|
	\leq2\sup_{t\in I}|F(x,t)|$. Taking the supremum over partitions and
	then applying H\"older's inequalit, we can prove (\refeq{eq12}). The last assertion follows
	from the corresponding inclusion for sequence norms.
\end{proof}

\begin{lemma}[Du--Zhang {\cite[Corollary 1.7]{DZ2019}}]
\label{le1}
Let $R\geq1$ and let $X$ be a union of lattice unit cubes in
$B^{n+1}(0,CR)$. Let
$$
 \gamma(X)=\sup_{\substack{z\in\R^{n+1}\\1\leq\rho\leq CR}}
 \frac{\#\{B:B\text{ is a unit cube of }X,\ B\subset B^{n+1}(z,\rho)\}}
 {\rho^n}.
$$
If $\operatorname{supp}\widehat f\subset B^n(0,2)$, then
\begin{equation}\label{eq7}
 \|e^{it\Delta}f\|_{L^2(X)}
 \lesssim_{\varepsilon,C}
 \gamma(X)^{1/(n+1)}R^{n/(2(n+1))+\varepsilon}\|f\|_2.
\end{equation}
\end{lemma}

\begin{lemma}{\cite[Theorems 1.1 and 1.2]{DL2025}}
\label{le2}
Let $R\geq1$ and let $\mu$ be a positive Borel measure satisfying
$$
 \operatorname{supp}\mu\subset B^{n+1}(0,CR),\qquad
 \mu(B^{n+1}(z,\rho))\leq C\rho^n
 \quad(z\in\R^{n+1},\ \rho>0).
$$
If $\operatorname{supp}\widehat g\subset B^n(0,2)$, then
\begin{equation}\label{eq8}
 \|e^{it\Delta}g\|_{L^{d_n}(d\mu)}
 \lesssim_{\varepsilon,C}R^\varepsilon\|g\|_2,
\end{equation}
where $d_n$ is defined in \eqref{eq2}.
For $d_n\leq p\leq\infty$, one also has
\begin{eqnarray}
 \|V^\infty_{[0,R]}\{e^{it\Delta}g\}\|_{L^p(B^n(0,R))}
 =2\left\|\sup_{0\leq t\leq R}|e^{it\Delta}g|\right\|_{L^p(B^n(0,R))}
 \nonumber\lesssim_{\varepsilon,p}R^\varepsilon\|g\|_2.
 \label{eq9}
\end{eqnarray}
\end{lemma}

For $N,R\geq1$ and $z_0=(x_0,t_0)$, set
$$
 \mathcal P_R^N(z_0)=
 \{(x,t):|t-t_0|\leq R,\ x-x_0-2N(t-t_0)e_1\in Q_R\},
 \qquad \mathcal P_R^N=\mathcal P_R^N(0),
$$
where $Q_R=[-R/2,R/2]^n$. 
The next lemma follows from the  
\cite[Proposition A.1]{LMY2024} and the bilinear-to-linear
criterion of \cite{LRV2011} which is stated in
\cite[Proposition 5.1]{LMY2024}. We give the interpolation
and rescaling argument explicitly.

\begin{lemma}\label{le3}
Let
$$
 \frac{2(n+1)}n<p_0<\frac{2(n+3)}{n+1},\qquad
 \frac{n+1}{p_0}+\frac1{q_0}=\frac n2.
$$
Suppose that, for every $\varepsilon>0$,
\begin{equation}\label{eq10}
 \big\||e^{it\Delta}f_1\,e^{it\Delta}f_2|^{1/2}\big\|
 _{L_x^{p_0}L_t^{q_0}(\mathcal P_R^N(z_0))}
 \lesssim_\varepsilon N^{1/p_0-1/q_0}R^\varepsilon
 \|f_1\|_2^{1/2}\|f_2\|_2^{1/2}
\end{equation}
uniformly in $N,R\geq1$ and $z_0\in\R^{n+1}$ whenever
$$
 \operatorname{supp}\widehat f_i\subset B^n(Ne_1,1),\qquad
 \operatorname{dist}(\operatorname{supp}\widehat f_1,
                     \operatorname{supp}\widehat f_2)\geq\tfrac12.
$$
Then, for
$$
 p_0<P<\frac{2(n+2)}n,\qquad
 \frac{n+1}{P}+\frac1Q=\frac n2,
$$
one has
\begin{equation}\label{eq11}
 \|e^{it\Delta}f\|_{L_x^PL_t^Q(\R^{n+1})}\lesssim \|f\|_2
\end{equation}
for $\operatorname{supp}\widehat f\subset\{1/2\leq|\xi|\leq2\}$.
\end{lemma}

\begin{proof}
	 For each critical
	pair $(P,Q)$ in the stated range, we also use the parabolic
	rescaling and the interpolation method in
	\cite[Section 5, (5.2)--(5.4)]{LMY2024},  which gives the bilinear estimate \cite[(5.5)]{LMY2024}
	throughout a neighborhood of $(P,Q)$. For $2<P<Q<\infty$, by \cite[Proposition 5.1]{LMY2024}, we can obtain
    $$\|e^{it\Delta}f\|_{L_x^P L_t^Q(\mathbb R^{n+1})}
	\lesssim \|f\|_{\dot H^{\,n/2-n/P-2/Q}}\lesssim \|f\|_{2}.$$
	
\end{proof}

We first apply the weighted restriction estimate of Du and Li
\cite[Theorem 1.2]{DL2025}, recalled in Lemma \ref{le2}.
A Galilean transformation and a weighted formulation of the time
supremum yield the following maximal estimate for initial data
with Fourier support in $B^n(Ne_1,1)$.
This improves the spatial exponent in
\cite[Theorem 1.4]{LMY2024} from
$2+4/(n+1+\sum_{j=2}^n j^{-1})$ to $d_n$,
while retaining the factor $N^{1/d_n}R^\varepsilon$.

\begin{lemma}\label{le8}
	Let $N,R\geq1$, and  $f\in\mathcal S(\mathbb R^n)$ satisfy
	$$
	\operatorname{supp}\mathcal F_x(f)\subset B^n(Ne_1,1).
	$$
	Then, for every $\varepsilon>0$,
	\begin{equation}\label{eq26}
		\|e^{it\Delta}f\|_{L_x^{d_n}L_t^\infty(\mathcal P_R^N(z_0))}
		\lesssim_\varepsilon N^{1/d_n}R^\varepsilon\|f\|_2
	\end{equation}
	uniformly in $z_0\in\R^{n+1}$.
\end{lemma}

\begin{proof}
	Let  $g(x)=e^{-iNx_1}f(x)$. Then we obtain
	$$|e^{it\Delta}f(x)|=|(e^{it\Delta}g)(x-2Nte_1)|,\qquad\operatorname{supp}\widehat g\subset B^n(0,1).
	$$
	Let $w\geq0$ satisfy
	$$
	\operatorname{supp}w\subset\mathcal P_R^N,\qquad
	\operatorname*{ess\,sup}_{x\in\R^n}\int_\R w(x,t)\,dt\leq1,
	\qquad d\mu(y,t)=N^{-1}w(y+2Nte_1,t)\,dy\,dt.
	$$
	Then $\operatorname{supp}\mu\subset B^{n+1}(0,CR)$ and
	$\mu(B(z,\rho))\lesssim\rho^n$. Lemma \ref{le2} gives
	$$
	\int_{\mathcal P_R^N}|e^{it\Delta}f(x)|^{d_n}w(x,t)\,dx\,dt
	=N\int|e^{it\Delta}g(y)|^{d_n}\,d\mu(y,t)
	\lesssim_\varepsilon NR^{\varepsilon d_n}\|f\|_2^{d_n}.
	$$
	Taking the supremum over $w$ proves \eqref{eq26} by
	mixed-norm duality.
\end{proof}

We next combine Lemma \ref{le8} with Tao's bilinear restriction
estimate \cite[Corollary 1.3]{Tao2003} and follow the interpolation
argument of Li, Miao and Yu \cite[Section 5]{LMY2024}.
The $\varepsilon$-removal and bilinear-to-linear steps are
collected in Lemma \ref{le3}.
This gives the following critical mixed-norm estimate with
$L^2$ initial data supported in the unit frequency annulus.
The improved maximal estimate in Lemma \ref{le8} lowers the
spatial threshold to $G_n$, thereby extending the range in
the unit-frequency version of
\cite[Theorem 1.6(ii)]{LMY2024}.

\begin{lemma}\label{th3}
	Let $G_n<P\leq q_n$, and define $Q(P)$ by
	$$
	\frac{n+1}{P}+\frac1{Q(P)}=\frac n2.
	$$
	Suppose that $f\in\mathcal S(\mathbb R^n)$ satisfies
	$$
	\operatorname{supp}\mathcal F_x(f)
	\subset\left\{\xi\in\mathbb R^n:\frac12\leq|\xi|\leq2\right\}.
	$$
	Then
	\begin{equation}\label{eq27}
		\|e^{it\Delta}f\|_{L_x^P L_t^{Q(P)}(\R^{n+1})}
		\lesssim_P\|f\|_2.
	\end{equation}
\end{lemma}

\begin{proof}
	Let
	$
	a_n=\frac{2(n+3)}{n+1}$ and $\mathcal P=\mathcal P_R^N(z_0),$ where
	$$
	\operatorname{supp}\widehat f_i\subset B^n(Ne_1,1),
	\qquad
	\operatorname{dist}
	\bigl(\operatorname{supp}\widehat f_1,
	\operatorname{supp}\widehat f_2\bigr)\geq\frac12.
	$$
	By H\"older's inequality and Lemma \ref{le8},
	\begin{eqnarray*}
		\||e^{it\Delta}f_1\,e^{it\Delta}f_2|^{1/2}\|_{L_x^{d_n}L_t^\infty(\mathcal P)}
		&\leq&
		\|e^{it\Delta}f_1\|_{L_x^{d_n}L_t^\infty(\mathcal P)}^{1/2}
		\|e^{it\Delta}f_2\|_{L_x^{d_n}L_t^\infty(\mathcal P)}^{1/2}
		\\
		&\lesssim&
		\bigl(N^{1/d_n}R^\varepsilon\|f_1\|_2\bigr)^{1/2}
		\bigl(N^{1/d_n}R^\varepsilon\|f_2\|_2\bigr)^{1/2}
		\\
		&=&
		N^{1/d_n}R^\varepsilon
		\|f_1\|_2^{1/2}\|f_2\|_2^{1/2}.
	\end{eqnarray*} 
	Let $g_i(x)=e^{-iNx_1}f_i(x),\ i=1,2.$ Then
	$$
	\widehat g_i(\xi)=\widehat f_i(\xi+Ne_1),\qquad
	\operatorname{supp}\widehat g_i\subset B^n(0,1),\qquad
	\|g_i\|_2=\|f_i\|_2,\qquad
	\operatorname{dist}
	\bigl(\operatorname{supp}\widehat g_1,
	\operatorname{supp}\widehat g_2\bigr)\geq\frac12.
	$$
	Moreover,
	$$
	e^{it\Delta}f_i(x)
	=
	e^{iNx_1-itN^2}
	(e^{it\Delta}g_i)(x-2Nte_1).
	$$
	Fix $\varepsilon>0$ and choose $b>a_n$ sufficiently close to
	$a_n$ that
	$$
	(n+1)\left(\frac1{a_n}-\frac1b\right)<\varepsilon.
	$$
	Since $b/2>(n+3)/(n+1)$, Tao's bilinear restriction estimate
	\cite[Corollary 1.3]{Tao2003} gives
	$$
	\|e^{it\Delta}g_1\,e^{it\Delta}g_2\|
	_{L^{b/2}(\mathbb R^{n+1})}
	\lesssim_b \|g_1\|_2\|g_2\|_2.
	$$
	By variable change $y=x-2Nte_1$,
	\begin{eqnarray*}
		\||e^{it\Delta}g_1\,e^{it\Delta}g_2|^{1/2}\|_{L^b(\mathbb R^{n+1})}
		&=&
		\left(
		\int_{\mathbb R}\int_{\mathbb R^n}
		|e^{it\Delta}g_1(y)\,e^{it\Delta}g_2(y)|^{b/2}
		\,dy\,dt
		\right)^{1/b}
		\\
		&=&
		\|e^{it\Delta}g_1\,e^{it\Delta}g_2\|
		_{L^{b/2}(\mathbb R^{n+1})}^{1/2}
		\\
		&\lesssim&
		\|f_1\|_2^{1/2}\|f_2\|_2^{1/2}.
	\end{eqnarray*}
	Writing $z_0=(x_0,t_0)$, we have
	$$
	|\mathcal P|
	=
	\int_{t_0-R}^{t_0+R}
	|x_0+2N(t-t_0)e_1+Q_R|\,dt
	=
	2R^{n+1}.
	$$
	Therefore, H\"older's inequality yields
	\begin{eqnarray*}
		\||e^{it\Delta}f_1\,e^{it\Delta}f_2|^{1/2}\|_{L^{a_n}(\mathcal P)}
		&\leq&
		|\mathcal P|^{1/a_n-1/b}\||e^{it\Delta}f_1\,e^{it\Delta}f_2|^{1/2}\|_{L^b(\mathcal P)}
		\\
		&\leq&
		(2R^{n+1})^{1/a_n-1/b}
		\|F\|_{L^b(\mathbb R^{n+1})}
		\\
		&\lesssim&
		R^{(n+1)(1/a_n-1/b)}
		\|f_1\|_2^{1/2}\|f_2\|_2^{1/2}
		\\
		&\lesssim&
		R^\varepsilon
		\|f_1\|_2^{1/2}\|f_2\|_2^{1/2}.
	\end{eqnarray*}
	Choose $\theta$ so that
	$$\frac{(n+1)(1-\theta)}{d_n}+\frac{(n+2)\theta}{a_n}=\frac {n}{2},
	$$
	which implies
	$$
	\theta=\frac{(n+3)\bigl(nd_n-2(n+1)\bigr)}{(n+1)\bigl((n+2)d_n-2(n+3)\bigr)},
	\qquad
	1-\theta
	=
	\frac{2d_n}
	{(n+1)\bigl((n+2)d_n-2(n+3)\bigr)}.
	$$
	Since $2(n+1)/n<d_n<a_n$, we have $0<\theta<1$.
	Furthermore,
	\begin{eqnarray*}
		\frac{1-\theta}{d_n}+\frac{\theta}{a_n}=\frac{n(n+1)d_n-2(n^2+2n-1)}{2(n+1)\bigl((n+2)d_n-2(n+3)\bigr)}=\frac1{G_n},
	\end{eqnarray*}
	and
	\begin{eqnarray*}
		\frac{\theta}{a_n}=\frac{nd_n-2(n+1)}{2\bigl((n+2)d_n-2(n+3)\bigr)}=\frac n2-\frac{n+1}{G_n}=\frac1{Q_n}.
	\end{eqnarray*}
	In particular,
	$$
	\frac{1-\theta}{d_n}
	=\frac1{G_n}-\frac1{Q_n}.
	$$
	Thus we can obtain 
	$$
	\|\mathbf1_{\mathcal P}(x,\cdot)|e^{it\Delta}f_1\,e^{it\Delta}f_2|^{1/2}\|_{L_t^{Q_n}}
	\leq
	\|\mathbf1_{\mathcal P}(x,\cdot)|e^{it\Delta}f_1\,e^{it\Delta}f_2|^{1/2}\|_{L_t^\infty}^{1-\theta}
	\|\mathbf1_{\mathcal P}(x,\cdot)|e^{it\Delta}f_1\,e^{it\Delta}f_2|^{1/2}\|_{L_t^{a_n}}^\theta.
	$$
	Following H\"older's inequality, we obtain
	\begin{eqnarray}
		\||e^{it\Delta}f_1\,e^{it\Delta}f_2|^{1/2}\|_{L_x^{G_n}L_t^{Q_n}(\mathcal P)}
		&\leq&
		\||e^{it\Delta}f_1\,e^{it\Delta}f_2|^{1/2}\|_{L_x^{d_n}L_t^\infty(\mathcal P)}^{1-\theta}
		\||e^{it\Delta}f_1\,e^{it\Delta}f_2|^{1/2}\|_{L^{a_n}(\mathcal P)}^\theta
		\nonumber\\
		&\lesssim&
		N^{(1-\theta)/d_n}R^\varepsilon
		\|f_1\|_2^{1/2}\|f_2\|_2^{1/2}
		\nonumber\\
		&=&
		N^{1/G_n-1/Q_n}R^\varepsilon
		\|f_1\|_2^{1/2}\|f_2\|_2^{1/2}.
		\label{eq28}
	\end{eqnarray}
	Due to $
	\frac1{G_n}
	=
	\frac{1-\theta}{d_n}+\frac{\theta}{a_n}$ with $0<\theta<1$, we have
	$$
	\frac{2(n+1)}n<d_n<G_n<a_n,
	\qquad
	\frac{n+1}{G_n}+\frac1{Q_n}=\frac n2.
	$$
	By Lemma \ref{le3} to \eqref{eq28} with
	$(p_0,q_0)=(G_n,Q_n)$ gives \eqref{eq27} for
	$$
	G_n<P<\frac{2(n+2)}n.
	$$
	For $\frac{2(n+2)}n\leq P\leq q_n,$ choose $\lambda\in[0,1]$ such that
	$$
	\frac1P
	=
	(1-\lambda)\frac{n}{2(n+2)}
	+\lambda\frac{n-1}{2(n+1)}.
	$$
	Then
	$$
	\frac1{Q(P)}
	=
	\frac n2-\frac{n+1}{P}
	=
	(1-\lambda)\frac{n}{2(n+2)}+\frac{\lambda}{2}.
	$$
	By \cite[Theorem 1.2, (5)]{KT1998}, and \eqref{eq24}, we  obtain
	$$
	\|e^{it\Delta}f\|_{L^{2(n+2)/n}(\mathbb R^{n+1})}
	\lesssim\|f\|_2,
	\qquad
	\|e^{it\Delta}f\|_{L_x^{q_n}L_t^2(\mathbb R^{n+1})}
	\lesssim\|f\|_2.
	$$
	Hence, by mixed-norm interpolation,
	\begin{eqnarray*}
		\|e^{it\Delta}f\|_{L_x^PL_t^{Q(P)}(\mathbb R^{n+1})}
		&\leq&
		\|e^{it\Delta}f\|_{L^{2(n+2)/n}(\mathbb R^{n+1})}^{1-\lambda}
		\|e^{it\Delta}f\|_{L_x^{q_n}L_t^2(\mathbb R^{n+1})}^{\lambda}
		\\
		&\lesssim&
		\|f\|_2.
	\end{eqnarray*}a
	This proves \eqref{eq27}.
\end{proof}

\section{Local variational estimates for Schr\"odinger equation}
\label{sec2}

We first prove the local estimates at $(p,r)=(2,q_n)$, for
$1\leq r\leq2$, and on the curve $(n+1)/p+1/r=n/2$.
We then interpolate these estimates on $B^n(0,R)\times[0,R]$.
The passage to the global Sobolev estimate is made in Section \ref{sec4}.

Let $Q_j$ be the spatial unit cube centered at $j\in\Z^n$, and
let $I_m=[m,m+1)$ for $m\in\Z$. Let $B_{j,m}=Q_j\times I_m$
and define
\begin{equation}\label{eq13}
 f_{j,m}:=\left(\int_{B_{j,m}}|e^{it\Delta}f(x)|^2\,dx\,dt\right)^{1/2}.
\end{equation}

\begin{lemma}\label{le5}
Let $1\leq p\leq\infty$, $1\leq r<\infty$, and $R\geq1$.
Suppose that $f\in\mathcal S(\mathbb R^n)$ satisfies
$$
 \operatorname{supp}\mathcal F_x(f)
 \subset\left\{\xi\in\mathbb R^n:\frac12\leq|\xi|\leq2\right\}.
$$
Then
\begin{eqnarray}
 \lefteqn{\|V^r_{[0,R]}\{e^{it\Delta}f\}\|_{L^p(B^n(0,R))}}\nonumber\\
 &\lesssim&
 \left(\sum_{Q_j\cap B^n(0,4R)\neq\varnothing}
 \left(\sum_{I_m\cap[-4R,4R]\neq\varnothing} f_{j,m}^r\right)^{p/r}
 \right)^{1/p}
 +O(R^{-1000})\|f\|_2.
 \label{eq14}
\end{eqnarray}
For $p=\infty$, the outer sum is replaced by a supremum.
\end{lemma}

\begin{proof}
	By the \cite[(3.2)]{CT2026}, we can obtain
\begin{equation}\label{eq15}
	V^r_{[0,R]}\{e^{it\Delta}f(x)\}
	\lesssim
	\left(
	\sum_{I_m\cap[0,R]\neq\varnothing}
	\left[
	\sup_{t\in I_m}
	\left(
	|e^{it\Delta}f(x)|
	+|\partial_t e^{it\Delta}f(x)|
	\right)
	\right]^r
	\right)^{1/r}.
\end{equation}
By \eqref{eq15} and $|Q_j|=1$, we obtain
\begin{eqnarray}\label{eq16}
	&&
	\|V^r_{[0,R]}\{e^{it\Delta}f\}\|_{L^p(B^n(0,R))}
	\nonumber\\
	&\lesssim&
	\left(
	\sum_{Q_j\cap B^n(0,R)\neq\varnothing}
	\int_{Q_j}
	\left[
	\sum_{I_m\cap[0,R]\neq\varnothing}
	\left(
	\sup_{t\in I_m}
	\bigl(
	|e^{it\Delta}f(x)|
	+|\partial_t e^{it\Delta}f(x)|
	\bigr)
	\right)^r
	\right]^{p/r}
	dx
	\right)^{1/p}
	\nonumber\\
	&\leq&
	\left(
	\sum_{Q_j\cap B^n(0,R)\neq\varnothing}
	\left(
	\sum_{I_m\cap[0,R]\neq\varnothing}
	\left[
	\sup_{(x,t)\in B_{j,m}}
	\bigl(
	|e^{it\Delta}f(x)|
	+|\partial_t e^{it\Delta}f(x)|
	\bigr)
	\right]^r
	\right)^{p/r}
	\right)^{1/p}.
\end{eqnarray}

Since the space--time Fourier transform of $e^{it\Delta}f$
is supported on the paraboloid $\tau=-|\xi|^2$, and
$\operatorname{supp}\widehat f
\subset\{\xi\in\mathbb R^n:1/2\leq|\xi|\leq2\}$, we have
$$
\operatorname{supp}\mathcal F_{x,t}(e^{it\Delta}f)
\subset S
:=
\left\{
(\xi,\tau)\in\mathbb R^n\times\mathbb R:
1/2\leq|\xi|\leq2,\ \tau=-|\xi|^2
\right\}.
$$
Let $\psi\in C_c^\infty(\mathbb R^{n+1})$ such that
$\psi\equiv1$ on an open neighborhood of $S$, and set
$\eta=\mathcal F_{x,t}^{-1}\psi$. This gives
\begin{eqnarray}\label{eq17}
	\sup_{(x,t)\in B_{j,m}}|e^{it\Delta}f(x)|
	&\leq&
	\sum_{(j',m')\in\mathbb Z^{n+1}}
	\sup_{(x,t)\in B_{j,m}}
	\left|
	\int_{B_{j',m'}}
	e^{is\Delta}f(w)\eta(x-w,t-s)\,dw\,ds
	\right|
	\nonumber\\
	&\leq&
	\sum_{(j',m')\in\mathbb Z^{n+1}}
	\sup_{(x,t)\in B_{j,m}}
	\|\eta(x-\cdot,t-\cdot)\|_{L^2(B_{j',m'})}
	f_{j',m'}.
\end{eqnarray}
Define
$$
\widetilde f_{j,m}
=
\begin{cases}
	f_{j,m},
	&
	B_{j,m}\cap
	\bigl(B^n(0,4R)\times[-4R,4R]\bigr)
	\neq\varnothing,
	\\[4pt]
	0,&\text{otherwise}.
\end{cases}
$$
By \eqref{eq17} and the rapid decay of $\eta$, for any
$M>0$, we obtain
\begin{eqnarray*}
	\sup_{(x,t)\in B_{j,m}}|e^{it\Delta}f(x)|
	&\lesssim&
	\sum_{(j',m')\in\mathbb Z^{n+1}}
	(1+|(j-j',m-m')|)^{-M}f_{j',m'}
	\\
	&=&
	\sum_{Q_{j'}\cap B^n(0,4R)\neq\varnothing}
	\sum_{I_{m'}\cap[-4R,4R]\neq\varnothing}
	(1+|(j-j',m-m')|)^{-M}\widetilde f_{j',m'}
	+A_{j,m},
\end{eqnarray*}
where $A_{j,m}$ is remainder term which is defined as 
$$
A_{j,m}=\sum_{(j',m')\in\mathbb Z^{n+1}}
(1+|(j-j',m-m')|)^{-M}
\bigl(f_{j',m'}-\widetilde f_{j',m'}\bigr).
$$

For the main term, we have
\begin{eqnarray*}
	&&
	\left(
	\sum_{Q_j\cap B^n(0,R)\neq\varnothing}
	\left(
	\sum_{I_m\cap[0,R]\neq\varnothing}
	\left[
	\sum_{Q_{j'}\cap B^n(0,4R)\neq\varnothing}
	\sum_{I_{m'}\cap[-4R,4R]\neq\varnothing}
	(1+|(j-j',m-m')|)^{-M}\widetilde f_{j',m'}
	\right]^r
	\right)^{p/r}
	\right)^{1/p}
	\\
	&=&
	\left(
	\sum_{Q_j\cap B^n(0,R)\neq\varnothing}
	\left(
	\sum_{I_m\cap[0,R]\neq\varnothing}
	\left[
	\sum_{(j',m')\in\mathbb Z^{n+1}}
	(1+|(j',m')|)^{-M}
	\widetilde f_{j-j',m-m'}
	\right]^r
	\right)^{p/r}
	\right)^{1/p}
	\\
	&\leq&
	\sum_{(j',m')\in\mathbb Z^{n+1}}(1+|(j',m')|)^{-M}
	\left(
	\sum_{Q_j\cap B^n(0,R)\neq\varnothing}
	\left(
	\sum_{I_m\cap[0,R]\neq\varnothing}
	|\widetilde f_{j-j',m-m'}|^r
	\right)^{p/r}
	\right)^{1/p}
	\\
	&\leq&
	\left(\sum_{(j',m')\in\mathbb Z^{n+1}}(1+|(j',m')|)^{-M}\right)
	\left(
	\sum_{Q_j\cap B^n(0,4R)\neq\varnothing}
	\left(
	\sum_{I_m\cap[-4R,4R]\neq\varnothing}|\widetilde f_{j,m}|^r
	\right)^{p/r}
	\right)^{1/p}
	\\
	&\lesssim&
	\left(
	\sum_{Q_j\cap B^n(0,4R)\neq\varnothing}
	\left(
	\sum_{I_m\cap[-4R,4R]\neq\varnothing}
	f_{j,m}^r
	\right)^{p/r}
	\right)^{1/p},
\end{eqnarray*}
provided that $M>n+1$.

For the remainder term $A_{j,m}$, suppose that
$Q_j\cap B^n(0,R)\neq\varnothing$ and
$I_m\cap[0,R]\neq\varnothing$. Then
$$
B_{j',m'}\cap
\bigl(B^n(0,4R)\times[-4R,4R]\bigr)
=\varnothing
\quad\Longrightarrow\quad
|(j-j',m-m')|\gtrsim R.
$$
Since $f_{j',m'}\leq\|f\|_2$, we obtain
\begin{eqnarray*}
	A_{j,m}
	&\leq&
	\sum_{\substack{(j',m')\in\mathbb Z^{n+1}\\|(j-j',m-m')|\gtrsim R}}
	(1+|(j-j',m-m')|)^{-M}f_{j',m'}
	\\
	&\leq&
	\|f\|_2
	\sum_{\substack{(j',m')\in\mathbb Z^{n+1}\\|(j-j',m-m')|\gtrsim R}}
	(1+|(j-j',m-m')|)^{-M}
	\\
	&\lesssim&
	R^{n+1-M}\|f\|_2.
\end{eqnarray*}
Note that 
$$
\begin{cases}
\#\{j:Q_j\cap B^n(0,R)\neq\varnothing\}\lesssim R^n
	\\[4pt]
	\#\{m:I_m\cap[0,R]\neq\varnothing\}\lesssim R,
\end{cases}
$$
Taking $M=1000+(n/p+1/r+n+1) $, we can deduce 
\begin{eqnarray*}
\left(\sum_{Q_j\cap B^n(0,R)\neq\varnothing}\left(\sum_{I_m\cap[0,R]\neq\varnothing}A_{j,m}^r\right)^{p/r}\right)^{1/p}\lesssim
	R^{n/p+1/r+n+1-M}\|f\|_2\lesssim 	R^{-1000}\|f\|_2.
\end{eqnarray*}

For the $\partial_t e^{it\Delta}f(x)$, we also have $\partial_t e^{it\Delta}f(x)=(\partial_t\eta)*e^{it\Delta}f(x),$
where the convolution is taken in space and time,
and $\partial_t\eta\in\mathcal S(\R^{n+1})$.
The same argument applies to $\partial_t e^{it\Delta}f(x)$.
Combining these estimates with \eqref{eq16} proves the lemma.
\end{proof}

\begin{lemma}\label{le6}
Let $q_n\leq r<\infty$, $R\geq1$, and let
$f\in\mathcal S(\mathbb R^n)$ satisfy
$$
 \operatorname{supp}\mathcal F_x(f)
 \subset\left\{\xi\in\mathbb R^n:\frac12\leq|\xi|\leq2\right\}.
$$
Then, for every $\varepsilon>0$,
\begin{equation}\label{eq18}
 \left(\sum_{Q_j\cap B^n(0,4R)\neq\varnothing}
 \left(\sum_{I_m\cap[-4R,4R]\neq\varnothing}f_{j,m}^{r}\right)^{2/r}
 \right)^{1/2}
 \lesssim_\varepsilon R^{n/(2(n+1))+\varepsilon}\|f\|_2.
\end{equation}
\end{lemma}

\begin{proof}
     Let
    $$
	\mathcal Q_R=\{Q_j:Q_j\cap B^n(0,4R)\neq\varnothing\},
	\qquad
	\mathcal I_R=\{m\in\mathbb Z:I_m\cap[-4R,4R]\neq\varnothing\}.
	$$
	Split $\mathcal I_R$ into $m\geq0$ and $m\leq-1$, and let
	$\mathcal M_R$ denote either part. Write $N_R=\#\mathcal M_R\lesssim R$.
	For dyadic $\lambda>0$ and dyadic integers $K\geq1$, put
$$
\left\{
\begin{aligned}
	&M(\lambda,Q_j)=\{m\in\mathcal M_R:f_{j,m}\geq\lambda\},\\
	&N(\lambda,Q_j)=\#M(\lambda,Q_j),\\
	&E_{\lambda,K}=\{Q_j\in\mathcal Q_R:K\leq N(\lambda,Q_j)<2K\},\\
	&X_{\lambda,K}=\bigcup_{Q_j\in E_{\lambda,K}}
	\bigcup_{m\in M(\lambda,Q_j)}B_{j,m}.
\end{aligned}
\right.
$$
	Since $r\geq q_n$, we have
	$$
	\left(\sum_{m\in\mathcal M_R} f_{j,m}^{r}\right)^{1/r}
	\leq
	\left(\sum_{m\in\mathcal M_R} f_{j,m}^{q_n}\right)^{1/q_n}.
	$$
	For each $Q_j\in\mathcal Q_R$, rearrange $\{f_{j,m}:m\in\mathcal M_R\}$
	in decreasing order as $f_{Q_j}(1)\geq\cdots\geq f_{Q_j}(N_R)$, and set
	\begin{equation}\label{eq20}
		F(Q_j)=
		\max_{\substack{1\leq K\leq N_R\\K\ \mathrm{dyadic}}}
		K^{1/q_n} f_{Q_j}(K).
	\end{equation}
	For $1\leq l\leq N_R$, choose a dyadic $K$ with $K\leq l<2K$. Then
	$$
	f_{Q_j}(l)\leq f_{Q_j}(K)
	\leq F(Q_j)K^{-1/q_n}
	\leq 2^{1/q_n}F(Q_j)l^{-1/q_n}.
	$$
	Summing over $l$ yields
	\begin{equation}\label{eq21}
		\left(\sum_{m\in\mathcal M_R}f_{j,m}^{q_n}\right)^{1/q_n}
		\leq 2^{1/q_n}F(Q_j)
		\left(\sum_{l=1}^{N_R}\frac{1}{l}\right)^{1/q_n}
		\lesssim (\log(2+R))^{1/q_n} F(Q_j).
	\end{equation}
	
	For $\mu>0$, let $Q_\mu=\{Q_j\in\mathcal Q_R:F(Q_j)\geq\mu\}$.
	Unitarity gives $f_{j,m}\leq\|f\|_2$, so
	$F(Q_j)\leq A_R:=N_R^{1/q_n}\|f\|_2$. Let
	$\eta=\min\{\varepsilon/2,1/4\}$ and take
	\[
	\mu_0=R^{-n^2/(2(n+1))+\eta}\|f\|_2.
	\]
	Since $\#\mathcal Q_R\lesssim R^n$, the layer-cake formula gives
	\begin{align*}
		\sum_{Q_j\in\mathcal Q_R}F(Q_j)^2
		=2\int_0^{A_R}\mu\,\#Q_\mu\,d\mu\lesssim R^n\mu_0^2+\int_{\mu_0}^{A_R}\mu\,\#Q_\mu\,d\mu .
	\end{align*}
	The first term is $R^{n/(n+1)+2\eta}\|f\|_2^2$.
	
	Fix $\mu_0<\mu\leq A_R$. For each dyadic $1\leq K\leq N_R$,
	choose $\lambda_K$ dyadic with
	\[
	\lambda_K\leq K^{-1/q_n}\mu<2\lambda_K .
	\]
	If $Q_j\in Q_\mu$, some such $K$ satisfies
	$f_{Q_j}(K)\geq K^{-1/q_n}\mu\geq\lambda_K$. Hence at least $K$
	entries in that column are at least $\lambda_K$. It follows that
	\[
	Q_\mu\subset
	\bigcup_{\substack{1\leq K\leq N_R\\K\ \mathrm{dyadic}}}
	\widetilde E_{\lambda_K,K},
	\qquad
	\widetilde E_{\lambda,K}=\{Q_j\in\mathcal Q_R:N(\lambda,Q_j)\geq K\}.
	\]
	
	For $Q_j\in E_{\lambda,K}$, the set $X_{\lambda,K}$ contains between
	$K$ and $2K$ unit cubes above $Q_j$. By disjointness of their interiors,
	\[
	\|e^{it\Delta}f\|_{L^2(X_{\lambda,K})}^2
	=\sum_{Q_j\in E_{\lambda,K}}\sum_{m\in M(\lambda,Q_j)}f_{j,m}^2
	\geq \lambda^2 K\,\#E_{\lambda,K}.
	\]
	The spatial projection of a ball of radius $\rho\geq1$ meets
	$O(\rho^n)$ spatial unit cubes. Thus a ball of radius $\rho$ contains
	$O(K\rho^n)$ unit cubes of $X_{\lambda,K}$, and
	$\gamma(X_{\lambda,K})\lesssim K$. Lemma \ref{le1}, with $\eta$
	in place of $\varepsilon$, gives
	\[
	\lambda K^{1/2}(\#E_{\lambda,K})^{1/2}
	\leq\|e^{it\Delta}f\|_{L^2(X_{\lambda,K})}
	\lesssim_\eta K^{1/(n+1)} R^{n/(2(n+1))+\eta}\|f\|_2 .
	\]
	Since $1/2-1/(n+1)=1/q_n$, this yields
	\begin{equation}\label{eq19}
		\lambda K^{1/q_n}(\#E_{\lambda,K})^{1/2}
		\lesssim_\eta R^{n/(2(n+1))+\eta}\|f\|_2 .
	\end{equation}
	
	For dyadic $K$, the sets $E_{\lambda,2^lK}$ with $l\geq0$
	partition $\widetilde E_{\lambda,K}$. By \eqref{eq19},
	\begin{align*}
		\#\widetilde E_{\lambda,K}
		&\lesssim_\eta
		\sum_{l\geq0}
		\left(\frac{R^{n/(2(n+1))+\eta}\|f\|_2}
		{\lambda (2^lK)^{1/q_n}}\right)^2\\
		&\lesssim_\eta
		\left(\frac{R^{n/(2(n+1))+\eta}\|f\|_2}
		{\lambda K^{1/q_n}}\right)^2 .
	\end{align*}
	The series converges because $q_n<\infty$. There are
	$O(\log(2+R))$ dyadic choices of $K$ and
	$\lambda_K K^{1/q_n}>\mu/2$. Consequently,
	\[
	\#Q_\mu\lesssim_\eta
	\log(2+R)
	\left(\frac{R^{n/(2(n+1))+\eta}\|f\|_2}{\mu}\right)^2 .
	\]
	Inserting this bound into the layer-cake integral gives
	\[
	\sum_{Q_j\in\mathcal Q_R} F(Q_j)^2
	\lesssim_\eta R^{n/(n+1)+2\eta}\|f\|_2^2
	\left(1+\log(2+R)\int_{\mu_0}^{A_R}\frac{d\mu}{\mu}\right).
	\]
	Here $A_R/\mu_0=N_R^{1/q_n}R^{n^2/(2(n+1))-\eta}$, so the integral is
	$O(\log(2+R))$. Taking square roots and using \eqref{eq21} proves
	the mixed-norm bound on $\mathcal M_R$, with a factor
	$(\log(2+R))^{1+1/q_n}$ multiplying
	$R^{n/(2(n+1))+\eta}\|f\|_2$. Applying this bound to the two time
	ranges and using the triangle inequality, we obtain
	\[
	\left(\sum_{Q_j\cap B^n(0,4R)\neq\varnothing}
	\left(\sum_{I_m\cap[-4R,4R]\neq\varnothing}
	f_{j,m}^{r}\right)^{2/r}\right)^{1/2}
	\lesssim_\varepsilon R^{n/(2(n+1))+\varepsilon}\|f\|_2,
	\]
	where the logarithmic factor is absorbed into
	$R^{\varepsilon-\eta}$.
\end{proof}

For $R\geq1$ and $f$ as in Lemma \ref{le6}, combining
Lemmas \ref{le5} and \ref{le6} at $r=q_n$ gives
\begin{equation}\label{eq22}
	\|V^{q_n}_{[0,R]}\{e^{it\Delta}f\}\|_{L^2(B^n(0,R))}
	\lesssim_\varepsilon R^{b_n/2+\varepsilon}\|f\|_2 .
\end{equation}
For $q_n\leq r<\infty$, monotonicity in $r$ gives
\[
\|V^r_{[0,R]}\{e^{it\Delta}f\}\|_{L^2(B^n(0,R))}
\leq\|V^{q_n}_{[0,R]}\{e^{it\Delta}f\}\|_{L^2(B^n(0,R))}
\lesssim_\varepsilon R^{b_n/2+\varepsilon}\|f\|_2 .
\]
For $r=\infty$, use
\[
\sup_{0\leq t\leq R}|e^{it\Delta}f(x)|
\leq |f(x)|+V^{q_n}_{[0,R]}\{e^{it\Delta}f\}(x).
\]
Together with \eqref{eq22}, this gives the maximal estimate at $p=2$.

\begin{lemma}\label{le7}
Let $2\leq p\leq\infty$, $1\leq r\leq2$, and $R\geq1$.
Suppose that $f\in\mathcal S(\mathbb R^n)$ satisfies
$$
 \operatorname{supp}\mathcal F_x(f)
 \subset\left\{\xi\in\mathbb R^n:\frac12\leq|\xi|\leq2\right\}.
$$
Then
\begin{equation}\label{eq23}
 \|V^r_{[0,R]}\{e^{it\Delta}f\}\|_{L^p(B^n(0,R))}
 \lesssim R^{1/r-1/2+\max\{\frac{n+1}{2p}-\frac{n-1}{4},\,0\}}\|f\|_2.
\end{equation}
\end{lemma}

\begin{proof}
	Since
	\[
	\#\{m\in\mathbb Z:I_m\cap[-4R,4R]\neq\varnothing\}
	\lesssim R,
	\]
	the finite-sequence inequality gives
	\[
	\left(
	\sum_{I_m\cap[-4R,4R]\neq\varnothing}f_{j,m}^r
	\right)^{1/r}
	\lesssim R^{1/r-1/2}
	\left(
	\sum_{I_m\cap[-4R,4R]\neq\varnothing}f_{j,m}^2
	\right)^{1/2}.
	\]
	By Vega's endpoint mixed-norm estimate \cite{V1988} and the
	Fourier support of $f$, we have
	\begin{equation}\label{eq24}
		\|e^{it\Delta}f\|_{L_x^{q_n}(\mathbb R^n;L_t^2(\mathbb R))}
		\lesssim
		\|f\|_{\dot H^{-1/(n+1)}(\mathbb R^n)}
		\lesssim \|f\|_2.
	\end{equation}
	Let $J_R=[-4R-1,4R+1]$. By 	Plancherel's theorem,
	\begin{align*}
		\|e^{it\Delta}f\|_{L_x^2(\mathbb R^n;L_t^2(J_R))}^2=
		|J_R|\|f\|_2^2
		=O(R)\|f\|_2^2
		\lesssim R\|f\|_2^2.
	\end{align*}
	
	Case 1: $2<p<q_n$. Choose $\theta\in(0,1)$ such that
    $$\frac1p=\frac{\theta}{2}+\frac{1-\theta}{q_n}.$$
	Since $q_n=2(n+1)/(n-1)$, this gives
	$$
	\theta
	=
	\frac{\frac1p-\frac1{q_n}}{\frac12-\frac1{q_n}}
	=\frac{n+1}{p}-\frac{n-1}{2},
	\frac{\theta}{2}=\frac{n+1}{2p}-\frac{n-1}{4}.
	$$
	By H\"older's inequality, we have 
    $$\frac{2}{\theta p}
	\quad\text{and}\quad
	\frac{q_n}{(1-\theta)p},
	\qquad
	\frac{\theta p}{2}+\frac{(1-\theta)p}{q_n}=1,$$
	which implies 
	\begin{eqnarray*}
		\|e^{it\Delta}f\|_{L_x^p(\mathbb R^n;L_t^2(J_R))}
		& \leq&
		\|e^{it\Delta}f\|_{L_x^2(\mathbb R^n;L_t^2(J_R))}^{\theta}
		\|e^{it\Delta}f\|_{L_x^{q_n}(\mathbb R^n;L_t^2(J_R))}^{1-\theta}
		\\
		&\leq&
		\|e^{it\Delta}f\|_{L_x^2(\mathbb R^n;L_t^2(J_R))}^{\theta}
		\|e^{it\Delta}f\|_{L_x^{q_n}(\mathbb R^n;L_t^2(\mathbb R))}^{1-\theta}
		\\
		& \lesssim&
		\bigl(R^{1/2}\|f\|_2\bigr)^\theta
		\|f\|_2^{1-\theta}
		=
		R^{\frac{n+1}{2p}-\frac{n-1}{4}}\|f\|_2.
	\end{eqnarray*}
    Then we obtain
	\begin{equation}\label{eq25}
		\|e^{it\Delta}f\|_{L_x^p(\mathbb R^n;L_t^2(J_R))}
		\lesssim
		R^{\frac{n+1}{2p}-\frac{n-1}{4}}\|f\|_2,
		\qquad 2\leq p\leq q_n.
	\end{equation}
	
    Case 2: $q_n\leq p\leq\infty$.
	Choose a fixed $K\in\mathcal S(\mathbb R^n)$ such that
	$\widehat K(\xi)=1$ for $|\xi|\leq2$. Since
	$$
	\mathcal F_x(e^{it\Delta}f)(\xi)
	=e^{-it|\xi|^2}\widehat f(\xi),
	\qquad
	\operatorname{supp}\widehat f
	\subset\{\xi\in\mathbb R^n:|\xi|\leq2\},
	$$
	we have, for every $t\in\mathbb R$,
	$$
	e^{it\Delta}f(x)
	=
	\int_{\mathbb R^n}K(x-y)e^{it\Delta}f(y)\,dy.
	$$
    By Minkowski's   inequality, we can deduce
	\begin{eqnarray*}
		 \left(
		\int_{J_R}|e^{it\Delta}f(x)|^2\,dt
		\right)^{1/2} 
		&=&
		\left(
		\int_{J_R}
		\left|
		\int_{\mathbb R^n}K(x-y)e^{it\Delta}f(y)\,dy
		\right|^2dt
		\right)^{1/2}
		\\
		&\leq&
		\int_{\mathbb R^n}|K(x-y)|
		\left(
		\int_{J_R}|e^{it\Delta}f(y)|^2\,dt
		\right)^{1/2}dy.
	\end{eqnarray*}
	Choose $\delta$ so that
	$$
	\frac1\delta=1+\frac1p-\frac1{q_n},
	\qquad
	1\leq\delta\leq\frac{q_n}{q_n-1}.$$
	Then we obtain
	inequality yields
	\begin{eqnarray*}
	 \left\|
		\left(
		\int_{J_R}|e^{it\Delta}f(x)|^2\,dt
		\right)^{1/2}
		\right\|_{L_x^p(\mathbb R^n)}
		&\leq&
		\left\|
		\int_{\mathbb R^n}|K(x-y)|
		\left(
		\int_{J_R}|e^{it\Delta}f(y)|^2\,dt
		\right)^{1/2}dy
		\right\|_{L_x^p(\mathbb R^n)}
		\\
		&\leq&
		\|K\|_{L^\delta(\mathbb R^n)}
		\left\|
		\left(
		\int_{J_R}|e^{it\Delta}f(x)|^2\,dt
		\right)^{1/2}
		\right\|_{L_x^{q_n}(\mathbb R^n)}\\
			& \lesssim&
		\left\|\left(\int_{J_R}|e^{it\Delta}f(x)|^2\,dt\right)^{1/2}
		\right\|_{L_x^{q_n}(\mathbb R^n)}.
	\end{eqnarray*}
    which indicates 
	\begin{eqnarray*}
		\|e^{it\Delta}f\|_{L_x^p(\mathbb R^n;L_t^2(J_R))}
		&\lesssim&
		\|e^{it\Delta}f\|_{L_x^{q_n}(\mathbb R^n;L_t^2(J_R))}
		\\
		&\leq&
		\|e^{it\Delta}f\|_{L_x^{q_n}(\mathbb R^n;L_t^2(\mathbb R))}
		\\
		&\lesssim& \|f\|_2
	\end{eqnarray*}
    for $q_n\leq p\leq\infty.$	Combining this with \eqref{eq25}, we obtain
	$$
	\|e^{it\Delta}f\|_{L_x^p(\mathbb R^n;L_t^2(J_R))}
	\lesssim
	R^{\max\{\frac{n+1}{2p}-\frac{n-1}{4},\,0\}}\|f\|_2,
	\qquad 2\leq p\leq\infty.$$
	Finally, since $|Q_j|=1$ and $p\geq2$, we have 
	$$
	\left(
	\sum_{Q_j\cap B^n(0,4R)\neq\varnothing}
	\left(
	\sum_{I_m\cap[-4R,4R]\neq\varnothing}f_{j,m}^2
	\right)^{p/2}
	\right)^{1/p}
	\leq
	\|e^{it\Delta}f\|_{L_x^p(\mathbb R^n;L_t^2(J_R))}.
	$$
     Combining these estimates with Lemma \ref{le5}, we conclude that
	\begin{eqnarray*}
		\|V^r_{[0,R]}\{e^{it\Delta}f\}\|_{L^p(B^n(0,R))}
		&\lesssim&
		R^{1/r-1/2}
		\|e^{it\Delta}f\|_{L_x^p(\mathbb R^n;L_t^2(J_R))}
		\\
		&\lesssim&
		R^{1/r-1/2+\max\{\frac{n+1}{2p}-\frac{n-1}{4},\,0\}}
		\|f\|_2.
	\end{eqnarray*}
	This proves \eqref{eq23}.
\end{proof}

\begin{lemma}\label{le9}
Let $G_n<p\leq q_n$, $R\geq1$, and let
$f\in\mathcal S(\mathbb R^n)$ satisfy
$$
 \operatorname{supp}\mathcal F_x(f)
 \subset\left\{\xi\in\mathbb R^n:\frac12\leq|\xi|\leq2\right\}.
$$
Let
$$
 q_p'=\frac{2p}{np-2(n+1)},\qquad
 \frac{n+1}{p}+\frac1{q_p'}=\frac n2.
$$
Then
\begin{equation}\label{eq29}
 \|V^{q_p'}_{[0,R]}\{e^{it\Delta}f\}\|_{L^p(B^n(0,R))}
 \lesssim_p\|f\|_2.
\end{equation}
\end{lemma}

\begin{proof}
Set
$$
\mathscr S=\{1\le j\le J:t_{j-1},t_j\in I_m
\text{ for some }m\in\mathbb Z\},
\qquad
\mathscr H=\{1,\ldots,J\}\setminus\mathscr S .
$$
Then we have 
$$
\sum_{j\in\mathscr S}|e^{it_{j}\Delta} f(x)-e^{it_{j-1}\Delta} f(x)|^{q_p'}
\le \sum_m\left(\int_{I_m}|\partial_{t}e^{it\Delta} f(x)|\,dt\right)^{q_p'}
\le \int_{\mathbb R}|\partial_{t}e^{it\Delta} f(x)|^{q_p'}\,dt .
$$
For $\mathscr H$, we obtain
\[
\sum_{j\in\mathscr H}|e^{it_{j}\Delta} f(x)-e^{it_{j-1}\Delta} f(x)|^{q_p'}
\lesssim_r \sum_m\sup_{t\in I_m}|e^{it\Delta} f(x)|^{q_p'}
\lesssim_r \int_{\mathbb R}|e^{it\Delta} f(x)|^r+|\partial_{t}e^{it\Delta} f(x)|^{q_p'}\,dt,
\]
which implies
\begin{eqnarray*}
	V_{[0,R]}^{q_p'}\{ e^{it\Delta} f(x)\}&\lesssim&\left(  \sum_{j\in\mathscr S}|e^{it_{j}\Delta} f(x)-e^{it_{j-1}\Delta} f(x)|^{q_p'}\right)^{1/q_p'} +\left( \sum_{j\in\mathscr H}|e^{it_{j}\Delta} f(x)-e^{it_{j-1}\Delta} f(x)|^{q_p'}\right)^{1/q_p'}\\
	&\lesssim&\| e^{it\Delta} f(x)\|_{L_t^{q_p'}(\mathbb R)}+\|\partial_{t}e^{it\Delta} f(x)\|_{L_t^{q_p'}(\mathbb R)}.
\end{eqnarray*}
By Lemma \ref{th3}, we can deduce 
\begin{eqnarray*}
	\|V^{q_p'}_{[0,R]}\{e^{it\Delta}f\}\|_{L^p(B^n(0,R))}\lesssim \| e^{it\Delta} f(x)\|_{L_{x}^{p}L_t^{q_p'}(\mathbb{R}^{n}\times\mathbb R)}+\|\partial_{t}e^{it\Delta} f(x)\|_{L_{x}^{p}L_t^{q_p'}(\mathbb{R}^{n}\times\mathbb R)}\lesssim \|f\|_{2}.
\end{eqnarray*}
\end{proof}

We use the interpolation inequality in Lemma \ref{le4}
on $B^n(0,R)\times[0,R]$ before passing to global estimates.

\begin{lemma}\label{le10}
Let $p>G_n$, $2\leq r\leq\infty$, and $R\geq1$.
Suppose that $f\in\mathcal S(\mathbb R^n)$ satisfies
$$
 \operatorname{supp}\mathcal F_x(f)
 \subset\left\{\xi\in\mathbb R^n:\frac12\leq|\xi|\leq2\right\}.
$$
Then, for every $\varepsilon>0$,
\begin{equation}\label{eq30}
 \|V^r_{[0,R]}\{e^{it\Delta}f\}\|_{L^p(B^n(0,R))}
 \lesssim_{p,r,\varepsilon}
 R^{\max\{0,\,(n+1)/(2p)+1/(2r)-n/4\}+\varepsilon}\|f\|_2.
\end{equation}
\end{lemma}

\begin{proof}
Suppose first that $G_n<p<q_n$ and use the exponent $q_p'$
from Lemma \ref{le9}.

\medskip\noindent
\textbf{Case 1: $2\leq r<q_p'$.}
Choose $0<\theta\leq1$ such that
$$
 \frac1r=\frac{\theta}{2}+\frac{1-\theta}{q_p'}.
$$
By Lemmas \ref{le7} and \ref{le9},
\begin{eqnarray*}
 \|V^r_{[0,R]}e^{it\Delta}f\|_{L^p(B^n(0,R))}
 &\leq&
 \|V^2_{[0,R]}e^{it\Delta}f\|_{L^p(B^n(0,R))}^{\theta}
 \|V^{q_p'}_{[0,R]}e^{it\Delta}f\|_{L^p(B^n(0,R))}^{1-\theta}\\
 &\lesssim& R^{\theta((n+1)/(2p)-(n-1)/4)}\|f\|_2.
\end{eqnarray*}
Since
$$
 \frac{n+1}{2p}-\frac{n-1}{4}
 =\frac12\left(\frac12-\frac1{q_p'}\right),
$$
the exponent is
$$
 \frac{\theta}{2}\left(\frac12-\frac1{q_p'}\right)
 =\frac12\left(\frac1r-\frac1{q_p'}\right)
 =\frac{n+1}{2p}+\frac1{2r}-\frac n4.
$$

\medskip\noindent
\textbf{Case 2: $q_p'\leq r<\infty$.}
By \eqref{eq29},
\begin{eqnarray*}
 \|V^r_{[0,R]}\{e^{it\Delta}f\}\|_{L^p(B^n(0,R))}
 \leq\|V^{q_p'}_{[0,R]}\{e^{it\Delta}f\}\|_{L^p(B^n(0,R))}\lesssim\|f\|_2.
\end{eqnarray*}
For $r=\infty$, Lemma \ref{le2} gives
\begin{eqnarray*}
	\|V^\infty_{[0,R]}\{e^{it\Delta}f\}\|_{L^p(B^n(0,R))}
	&\lesssim_{\varepsilon,p}&R^\varepsilon\|f\|_2
	=R^{\max\{0,\,(n+1)/(2p)-n/4\}+\varepsilon}\|f\|_2,
\end{eqnarray*}
since $p>G_n>2(n+1)/n$.
For $p\geq q_n$ and $2\leq r<\infty$, \eqref{eq23} gives
\begin{eqnarray*}
 \|V^r_{[0,R]}\{e^{it\Delta}f\}\|_{L^p(B^n(0,R))}
 \leq\|V^2_{[0,R]}\{e^{it\Delta}f\}\|_{L^p(B^n(0,R))}\lesssim&\|f\|_2.
\end{eqnarray*}
The case $r=\infty$ follows from Lemma \ref{le2}.
\end{proof}

\begin{lemma}\label{le11}
Let $2\leq p\leq G_n$, $R\geq1$, and let
$f\in\mathcal S(\mathbb R^n)$ satisfy
$$
 \operatorname{supp}\mathcal F_x(f)
 \subset\left\{\xi\in\mathbb R^n:\frac12\leq|\xi|\leq2\right\}.
$$
Let
$$
 \beta_n(p)=\frac{1-\theta_p}{2}b_n
 =\frac12\left(\frac1p+\frac1{\rho_n(p)}-A_n(p)\right).
$$
Then, for $2\leq r<\infty$ and every $\varepsilon>0$,
\begin{equation}\label{eq31}
 \|V^r_{[0,R]}\{e^{it\Delta}f\}\|_{L^p(B^n(0,R))}
 \lesssim_\varepsilon
 R^{\max\{\beta_n(p),\,(n+1)/(2p)+1/(2r)-n/4\}+\varepsilon}\|f\|_2.
\end{equation}
In particular,  for $2\leq r\leq\rho_n(p)$ the exponent is
$(n+1)/(2p)+1/(2r)-n/4$. The same formula holds at $p=G_n$
with $\beta_n(G_n)=0$.
\end{lemma}

\begin{proof}
Let $2<p\leq G_n$ and choose $P>G_n$ close to $G_n$.
Interpolate \eqref{eq22} with \eqref{eq29} at
$(P,Q(P))$. Choose $\theta$ so that
$$
 \frac1p=\frac{1-\theta}{2}+\frac{\theta}{P},\qquad
 \frac1{r_P}=\frac{1-\theta}{q_n}+\frac{\theta}{Q(P)}.
$$
By Lemma \ref{le4}, \eqref{eq22} and \eqref{eq29},
\begin{equation}\label{eq32}
 \|V^{r_P}_{[0,R]}\{e^{it\Delta}f\}\|_{L^p(B^n(0,R))}
 \lesssim_{\varepsilon,P}
 R^{(1-\theta)b_n/2+\varepsilon}\|f\|_2.
\end{equation}
The local exponent is $(1-\theta)b_n/2$. For  $P\downarrow G_n$, we have 
$$ r_P\longrightarrow\rho_n(p),\qquad
 (1-\theta)b_n/2\longrightarrow\beta_n(p).
$$
For $\beta_n(p)$ follows by interpolating
$$
 b_n=\frac12+\frac1{q_n},\qquad
 A_n(G_n)=\frac1{G_n}+\frac1{Q_n},
$$
and using that $A_n(p)$ is affine in $1/p$. Since $P>G_n$,
$$\theta=\frac{1/2-1/p}{1/2-1/P}<\frac{1/2-1/p}{1/2-1/G_n}=\theta_p.$$
The definitions of $r_P$ and $\rho_n(p)$ give
$$\frac1{r_P}-\frac1{\rho_n(p)}=b_n(\theta_p-\theta)>0,$$
so $r_P<\rho_n(p)$.
Fix $\varepsilon>0$ and choose $P>G_n$ sufficiently close
to $G_n$ that
$$
0<
\frac{(1-\theta)b_n}{2}-\beta_n(p)
=
\frac{b_n}{2}(\theta_p-\theta)
<
\frac{\varepsilon}{2}.
$$
By \eqref{eq32} with $\varepsilon/2$,
\begin{eqnarray*}
	\|V^{\rho_n(p)}_{[0,R]}e^{it\Delta}f\|_{L^p(B^n(0,R))}
	&\leq&
	\|V^{r_P}_{[0,R]}e^{it\Delta}f\|_{L^p(B^n(0,R))}
	\\
	&\lesssim&
	R^{(1-\theta)b_n/2+\varepsilon/2}\|f\|_2
	\\
	&\leq&
	R^{\beta_n(p)+\varepsilon}\|f\|_2.
\end{eqnarray*}
For $2\leq r<\rho_n(p)$, interpolate this estimate with
Lemma \ref{le7} at $r=2$. Then we have 
$$\frac{n+1}{2p}+\frac1{2r}-\frac {n}{4}.$$
For $\rho_n(p)<r<\infty$,
\begin{eqnarray*}
 \|V^r_{[0,R]}\{e^{it\Delta}f\}\|_{L^p(B^n(0,R))}
 &\leq&\|V^{\rho_n(p)}_{[0,R]}\{e^{it\Delta}f\}\|_{L^p(B^n(0,R))} \lesssim R^{\beta_n(p)+\varepsilon}\|f\|_2.
\end{eqnarray*}

For $p=2$ and $2\leq r\leq q_n$, choose $0\leq\vartheta\leq1$ so that
$$
 \frac1r=\frac{1-\vartheta}{2}+\frac{\vartheta}{q_n},
 \qquad
 \vartheta=\frac{1/2-1/r}{1/2-1/q_n}.
$$
Lemma \ref{le4}, \eqref{eq22} and \eqref{eq23} give
\begin{eqnarray*}
 \|V^r_{[0,R]}\{e^{it\Delta}f\}\|_{L^2(B^n(0,R))}
 &\leq&
 \|V^2_{[0,R]}\{e^{it\Delta}f\}\|_{L^2(B^n(0,R))}^{1-\vartheta}
 \|V^{q_n}_{[0,R]}\{e^{it\Delta}f\}\|_{L^2(B^n(0,R))}^{\vartheta}\\
 &\lesssim&
 R^{(1-\vartheta)/2+\vartheta b_n/2+\varepsilon}\|f\|_2\\
 &=&R^{1/4+1/(2r)+\varepsilon}\|f\|_2.
\end{eqnarray*}
For $p=2$ and $q_n<r<\infty$, \eqref{eq22} gives
\begin{eqnarray*}
 \|V^r_{[0,R]}\{e^{it\Delta}f\}\|_{L^2(B^n(0,R))}
 &\leq&\|V^{q_n}_{[0,R]}\{e^{it\Delta}f\}\|_{L^2(B^n(0,R))}\lesssim_\varepsilon R^{b_n/2+\varepsilon}\|f\|_2.
\end{eqnarray*}
Since
$$
 \rho_n(2)=q_n,\qquad
 \frac14+\frac1{2q_n}=\frac{b_n}{2}=\beta_n(2),
$$
these estimates prove \eqref{eq31}.
\end{proof}
For $d_n\leq p\leq G_n$, define $\eta_n(p)$ by
$$
\frac{1}{\eta_n(p)}
=
\frac{1/d_n-1/p}{1/d_n-1/G_n}\frac{1}{Q_n},
$$
with $\eta_n(d_n)=\infty$. In particular,
$\eta_n(G_n)=Q_n$.
\begin{lemma}\label{le12}
Let $R\geq1$, and let $f\in\mathcal S(\mathbb R^n)$ satisfy
$$
 \operatorname{supp}\mathcal F_x(f)
 \subset\left\{\xi\in\mathbb R^n:\frac12\leq|\xi|\leq2\right\}.
$$
Then, for every $\varepsilon>0$, the following estimates hold:
\begin{enumerate}[(i)]
\item If $d_n\leq p\leq G_n$ and $\eta_n(p)\leq r\leq\infty$, then
$$
 \|V^r_{[0,R]}\{e^{it\Delta}f\}\|_{L^p(B^n(0,R))}
 \lesssim_\varepsilon R^\varepsilon\|f\|_2.
$$
\item If $p=2$ and $q_n\leq r\leq\infty$, then
$$
 \|V^r_{[0,R]}\{e^{it\Delta}f\}\|_{L^2(B^n(0,R))}
 \lesssim_\varepsilon R^{b_n/2+\varepsilon}\|f\|_2.
$$
\end{enumerate}
\end{lemma}
\begin{proof}
	For (i), Lemma \ref{le2} gives
	$$
	\|V^\infty_{[0,R]}\{e^{it\Delta}f\}\|_{L^p(B^n(0,R))}
	\lesssim_\varepsilon R^\varepsilon\|f\|_2,
	\qquad p\geq d_n.
	$$
	This also proves the assertion when $p=d_n$, since
	$\eta_n(d_n)=\infty$. By
	$
	\rho_n(G_n)=Q_n,$ $
	\beta_n(G_n)=0,$ and $
	\frac{n+1}{2G_n}+\frac{1}{2Q_n}-\frac n4=0,
	$
	Lemma \ref{le11} yields
	$$
	\|V^{Q_n}_{[0,R]}\{e^{it\Delta}f\}\|_{L^{G_n}(B^n(0,R))}
	\lesssim_\varepsilon R^\varepsilon\|f\|_2.
	$$
    Let 
	$$
	\vartheta
	=
	\frac{1/d_n-1/p}{1/d_n-1/G_n}
	$$
	for $d_n<p\leq G_n$.	which indicates 
	$
	0<\vartheta\leq1,$ $
	\frac1p=\frac{1-\vartheta}{d_n}
	+\frac{\vartheta}{G_n},$
 and $\frac1{\eta_n(p)}=\frac{\vartheta}{Q_n}.$
	The variation interpolation inequality gives
	$$
	V^{\eta_n(p)}_{[0,R]}\{e^{it\Delta}f(x)\}
	\leq
	\bigl(V^\infty_{[0,R]}\{e^{it\Delta}f(x)\}\bigr)^{1-\vartheta}
	\bigl(V^{Q_n}_{[0,R]}\{e^{it\Delta}f(x)\}\bigr)^\vartheta.
	$$
	By H\"older's inequality, for
	$\eta_n(p)\leq r<\infty$,
	\begin{eqnarray*}
		\|V^r_{[0,R]}\{e^{it\Delta}f\}\|_{L^p(B^n(0,R))}
		&\leq&
		\|V^{\eta_n(p)}_{[0,R]}\{e^{it\Delta}f\}\|_{L^p(B^n(0,R))}
		\\
		&\leq&
		\|V^\infty_{[0,R]}\{e^{it\Delta}f\}\|_{L^{d_n}(B^n(0,R))}^{1-\vartheta}
		\|V^{Q_n}_{[0,R]}\{e^{it\Delta}f\}\|_{L^{G_n}(B^n(0,R))}^{\vartheta}
		\\
		&\lesssim&
		R^{(1-\vartheta)\varepsilon+\vartheta\varepsilon}
		\|f\|_2
		\\
		&=&
		R^\varepsilon\|f\|_2.
	\end{eqnarray*}
	This proves (i).
	
	For (ii), \eqref{eq22} gives, for $q_n\leq r<\infty$,
	\begin{eqnarray*}
		\|V^r_{[0,R]}\{e^{it\Delta}f\}\|_{L^2(B^n(0,R))}
		&\leq&
		\|V^{q_n}_{[0,R]}\{e^{it\Delta}f\}\|_{L^2(B^n(0,R))}
		\\
		&\lesssim&
		R^{b_n/2+\varepsilon}\|f\|_2.
	\end{eqnarray*}
	For $r=\infty$, we have 
	$$
	\sup_{0\leq t\leq R}|e^{it\Delta}f(x)|
	\leq |f(x)|+V^{q_n}_{[0,R]}\{e^{it\Delta}f(x)\}.
	$$
	Hence
	\begin{eqnarray*}
		\|V^\infty_{[0,R]}\{e^{it\Delta}f\}\|_{L^2(B^n(0,R))}
		&\leq&
		2\|f\|_{L^2(B^n(0,R))}
		+2\|V^{q_n}_{[0,R]}\{e^{it\Delta}f\}\|_{L^2(B^n(0,R))}
		\\
		&\lesssim&
		R^{b_n/2+\varepsilon}\|f\|_2.
	\end{eqnarray*}
	This proves (ii).
\end{proof}

\section{The proof of the sufficiency parts in Theorems \ref{th1} and \ref{th2}}
\label{sec4}

In this section, we first prove the local-to-global principle
for variation estimates. We then apply it to the local results
of Section \ref{sec2}. The spatial localization uses the
near--far decomposition in \cite{LRS2013,LRV2011,R2008};
we give the argument for the variation norm.

\begin{theorem}\label{th4}
Let $2\leq p\leq\infty$, $1\leq r\leq\infty$, and $\beta\geq0$.
Suppose that, for every $\varepsilon>0$, every $R\geq1$, and
every $g\in\mathcal S(\R^n)$ with
$$
 \operatorname{supp}\mathcal F_x(g)
 \subset\{\xi\in\R^n:1/2\leq|\xi|\leq2\},
$$
one has
\begin{equation}\label{eq33}
 \|V^r_{[0,R]}\{e^{it\Delta}g\}\|_{L^p(B^n(0,R))}
 \lesssim_\varepsilon R^{\beta+\varepsilon}\|g\|_2.
\end{equation}
Then, for every $s>A_n(p)+2\beta$,
$$
 \|V^r_{[0,1]}\{e^{it\Delta}f\}\|_{L^p(\R^n)}
 \lesssim_s\|f\|_{H^s(\R^n)}.
$$
\end{theorem}

\begin{proof}
We first localize in frequency and apply parabolic rescaling.
Let $g_R$ satisfy
$\operatorname{supp}\widehat g_R\subset\{R/2\leq|\xi|\leq2R\}$,
and define
$$
 g(x)=R^{-n/2}g_R(x/R).
$$
Then $\|g\|_2=\|g_R\|_2$, $\operatorname{supp}\mathcal F_x(g)
\subset\{\xi\in\R^n:1/2\leq|\xi|\leq2\}$, and
$$
 e^{it\Delta}g_R(x)=R^{n/2}(e^{iR^2t\Delta}g)(Rx).
$$
By variable change,
\begin{eqnarray*}
 \|V^r_{[0,1]}\{e^{it\Delta}g_R\}\|_{L^p(B^n(0,R))}
 &=&R^{n/2-n/p}
   \|V^r_{[0,R^2]}\{e^{it\Delta}g\}\|_{L^p(B^n(0,R^2))}\\
 &\lesssim&
 R^{A_n(p)+2\beta+\varepsilon}\|g_R\|_2,
\end{eqnarray*}
by \eqref{eq33} at radius $R^2$.
The same bound holds on translated balls and fixed enlargements
of the frequency annulus.

Choose $a\in C_c^\infty(\R^n)$ with $a=1$ on
$\{1/2\leq|\xi|\leq2\}$ and
$\operatorname{supp}a\subset\{1/4\leq|\xi|\leq4\}$.
Define $T_a$ by
$$
\widehat{T_a f}(\xi)
=a\left(\xi/R\right)\widehat f(\xi).$$
Partition $\R^n$ into cubes
$$
 Q_j=R(j+[0,1)^n),\qquad j\in\Z^n,
$$
and choose a smooth partition of unity $\{\chi_j\}_j$
adapted to these cubes, with uniformly bounded overlap.
Define
$$
 g_{R,j}:=T_a(\chi_jg_R)
 =\mathcal F_x^{-1}\bigl(a(\xi/R)\mathcal F_x(\chi_jg_R)\bigr).
$$
Since $T_ag_R=g_R$,
$$
 \sum_{j\in\mathbb Z^n}g_{R,j}=g_R,\qquad
 \operatorname{supp}\widehat g_{R,j}
 \subset\{\xi:R/4\leq|\xi|\leq4R\}.
$$
Plancherel and bounded overlap give
\begin{equation}\label{eq34}
 \sum_{j\in\mathbb Z^n}\|g_{R,j}\|_2^2
 \lesssim\sum_{j\in\mathbb Z^n}\|\chi_jg_R\|_2^2
 \lesssim\|g_R\|_2^2.
\end{equation}

Fix a sufficiently large dimensional constant $C_0$.
For $x\in Q_\ell$, decompose
$$
 g_R=\sum_{\|j-\ell\|_\infty\leq C_0}g_{R,j}
      +\sum_{\|j-\ell\|_\infty>C_0}g_{R,j}.
$$
For the near terms, cover $Q_\ell$ by a bounded number of
balls of radius $R$ and apply the translated local estimate.
Since $p\geq2$, bounded overlap and \eqref{eq34} give
\begin{eqnarray*}
 \left(\sum_{\ell\in\mathbb Z^n}
 \left\|V^r_{[0,1]}\left\{e^{it\Delta}
       \sum_{\|j-\ell\|_\infty\leq C_0}g_{R,j}\right\}
       \right\|_{L^p(Q_\ell)}^p\right)^{1/p} 
 &\lesssim_\varepsilon&
 R^{A_n(p)+2\beta+\varepsilon}
 \left(\sum_{\ell\in\mathbb Z^n}
  \left\|\sum_{\|j-\ell\|_\infty\leq C_0}g_{R,j}\right\|_2^2
       \right)^{1/2}\\
 &\lesssim_\varepsilon& R^{A_n(p)+2\beta+\varepsilon}\|g_R\|_2.
\end{eqnarray*}
 For the far terms and $\nu=0,1$, we have 
$$\partial_t^\nu e^{it\Delta}g_{R,j}(x)=\frac1{(2\pi)^n}\int_{\R^n}\int_{\R^n}(-i|\xi|^2)^\nu e^{i((x-y)\cdot\xi-t|\xi|^2)}a(\xi/R)\chi_j(y)g_R(y)\,dy\,d\xi.$$
If $x\in Q_\ell$, $y\in\operatorname{supp}\chi_j$,
$\|j-\ell\|_\infty>C_0$, and $0\leq t\leq1$, the phase is
nonstationary. By variable change $\xi=R\zeta$ and integration
by parts, for every $M$,
$$
 \sup_{\substack{x\in Q_\ell\\0\leq t\leq1}}
 |\partial_t^\nu e^{it\Delta}g_{R,j}(x)|
 \lesssim_{M,\nu}
 R^{3n/2+2\nu-2M}
 (1+\|j-\ell\|_\infty)^{-M}\|\chi_jg_R\|_2.
$$
Here we used $\|\chi_jg_R\|_1\lesssim R^{n/2}\|\chi_jg_R\|_2$.
For finite $r$,
$$
 V^r_{[0,1]}e^{it\Delta}g_{R,j}(x)
 \leq\int_0^1|\partial_te^{it\Delta}g_{R,j}(x)|\,dt.
$$
For $r=\infty$, add $2|g_{R,j}(x)|$ and use the same derivative
bound. Thus, in both cases,
$$\|V^r_{[0,1]}\{e^{it\Delta}g_{R,j}\}\|_{L^p(Q_\ell)}\lesssim_M R^{3n/2+2+n/p-2M}(1+\|j-\ell\|_\infty)^{-M}\|\chi_jg_R\|_2.$$
The discrete Young inequality, followed by $\ell^2\subset\ell^p$, bounds the sum of the far terms by $O(R^{-100})\|g_R\|_2$ after
choosing $M$ sufficiently large. Combining the two parts yields
\begin{equation}\label{eq35}
 \|V^r_{[0,1]}\{e^{it\Delta}g_R\}\|_{L^p(\R^n)}
 \lesssim_\varepsilon
 R^{A_n(p)+2\beta+\varepsilon}\|g_R\|_2.
\end{equation}

It remains to sum over frequencies. Write
$f=g_0+\sum_{k\geq1}g_k$ using a smooth inhomogeneous
Littlewood--Paley decomposition, with
$\operatorname{supp}\widehat g_k\subset\{|\xi|\approx2^k\}$.
Let $s_0=A_n(p)+2\beta$. For $s>s_0$, choose
$0<\varepsilon<s-s_0$. By subadditivity and Cauchy--Schwarz,
\begin{eqnarray*}
 \|V^r_{[0,1]}e^{it\Delta}f\|_p
 &\lesssim& \|g_0\|_2+
       \sum_{k\geq1}2^{k(s_0+\varepsilon)}\|g_k\|_2\\
 &\lesssim &
 \left(\|g_0\|_2^2+\sum_{k\geq1}2^{2ks}\|g_k\|_2^2\right)^{1/2}
 \lesssim_s\|f\|_{H^s}.
\end{eqnarray*}
For the low-frequency term, use
$V^r_{[0,1]}e^{it\Delta}g_0\leq\int_0^1|\Delta e^{it\Delta}g_0|\,dt$
and Bernstein's inequality; for $r=\infty$ also include
$2|g_0|$. Approximation by Schwartz functions completes the proof.
\end{proof}

\begin{theorem}\label{th5}
For each pair $(p,r)$ in Theorem \ref{th1} and every
$s>s(p,r)$,
$$
 \|V^r_{[0,1]}\{e^{it\Delta}f\}\|_{L^p(\R^n)}
 \lesssim_{p,r,s}\|f\|_{H^s(\R^n)}.
$$
\end{theorem}

\begin{proof}
We combine the local estimates of Section \ref{sec2}
with Theorem \ref{th4}. A local power $R^{\beta+\varepsilon}$
gives the global threshold
$$
 s_0=A_n(p)+2\beta.
$$

\medskip\noindent
\textbf{Case 1: $1\leq r\leq2$.}
Lemma \ref{le7} gives
$\beta=1/r-1/2+\max\{(n+1)/(2p)-(n-1)/4,0\}$.
Consequently,
$$
 s_0=\max\left\{\frac2r+\frac1p-\frac12,\,
                  A_n(p)+\frac2r-1\right\},
$$
which proves part (i).

\medskip\noindent
\textbf{Case 2: $p\geq G_n$ and $2\leq r\leq\infty$.}
For $p>G_n$, use Lemma \ref{le10}; for $p=G_n$
and finite $r$, use Lemma \ref{le11}.
The maximal boundary follows from Lemma \ref{le2}.
In each case we may take
$$
 \beta=\max\left\{0,\frac{n+1}{2p}+\frac1{2r}-\frac n4\right\}.
$$
Therefore
$$
 s_0=\max\left\{A_n(p),\frac1p+\frac1r\right\},
$$
proving part (ii).

\medskip\noindent
\textbf{Case 3: $2\leq p\leq G_n$ and $2\leq r\leq\rho_n(p)$.}
Lemma \ref{le11} gives
$$
 \beta=\frac{n+1}{2p}+\frac1{2r}-\frac n4.
$$
Thus $s_0=1/p+1/r$, as asserted in part (iii).

\medskip\noindent
\textbf{Case 4: $d_n\leq p\leq G_n$ and $\eta_n(p)\leq r\leq\infty$.}
Part (i) of Lemma \ref{le12} gives $\beta=0$.
Hence $s_0=A_n(p)$, proving part (iv).

\medskip\noindent
\textbf{Case 5: $p=2$ and $q_n\leq r\leq\infty$.}
Part (ii) of Lemma \ref{le12} gives $\beta=b_n/2$.
Since $A_n(2)=0$, we obtain $s_0=b_n$, proving part (v).
This completes the proof of sufficiency.
\end{proof}

\begin{theorem}
	Let $n\geq3$, $2\leq p\leq G_n$, and
	$\rho_n(p)\leq r\leq\infty$. Let
	\begin{equation}
		\sigma_n(p,r)=A_n(p)+b_n\max\left\{
		0,\,
		\frac{1/p-1/d_n}{1/2-1/d_n},\,
		\frac{1/p-1/d_n+(Q_n/r)(1/d_n-1/G_n)}
		{1/2-1/d_n+(Q_n/q_n)(1/d_n-1/G_n)}
		\right\}.
	\end{equation}
	If $s>\sigma_n(p,r)$, then
	$$
	\left\|V^r_{[0,1]}\{e^{it\Delta}f\}\right\|_{L^p(\R^n)}
	\lesssim_{p,r,s}\|f\|_{H^s(\R^n)}.
	$$
	
\end{theorem}

\begin{proof}
Let $R\geq1$ and
$\operatorname{supp}\widehat f\subset\{1/2\leq|\xi|\leq2\}$.
Lemma \ref{le11} at $(p,r)=(G_n,Q_n)$ and
Lemma \ref{le12} give
\begin{equation}\label{eq36}
 \left\{
 \begin{array}{rcl}
 \displaystyle
 \|V^{r_0}_{[0,R]}\{e^{it\Delta}f\}\|_{L^2(B^n(0,R))}
 &\lesssim_\varepsilon&
 \displaystyle R^{b_n/2+\varepsilon}\|f\|_2,
 \quad q_n\leq r_0\leq\infty,\\[6pt]
 \displaystyle
 \|V^{Q_n}_{[0,R]}\{e^{it\Delta}f\}\|_{L^{G_n}(B^n(0,R))}
 &\lesssim_\varepsilon&
 \displaystyle R^\varepsilon\|f\|_2,\\[6pt]
 \displaystyle
 \|V^\infty_{[0,R]}\{e^{it\Delta}f\}\|_{L^{d_n}(B^n(0,R))}
 &\lesssim_\varepsilon&
 \displaystyle R^\varepsilon\|f\|_2.
 \end{array}
 \right.
\end{equation}

\medskip\noindent
\textbf{Case 1: $2\leq p<d_n$ and
$1/r\leq(1/p-1/d_n)/(q_n(1/2-1/d_n))$.}
Let
$$
 \vartheta=\frac{1/p-1/d_n}{1/2-1/d_n},
 \qquad
 \frac1p=\frac{\vartheta}{2}+\frac{1-\vartheta}{d_n}.
$$
For finite $r$, set $r_0=\vartheta r$; for $r=\infty$,
set $r_0=\infty$. Then
$$
 0<\vartheta\leq1,\qquad
 r_0\geq q_n,\qquad
 \frac1r=\frac{\vartheta}{r_0}.
$$
Lemma \ref{le4} and \eqref{eq36} yield
\begin{eqnarray*}
 \|V^r_{[0,R]}\{e^{it\Delta}f\}\|_{L^p(B^n(0,R))}
 &\leq&
 \|V^{r_0}_{[0,R]}\{e^{it\Delta}f\}\|_{L^2(B^n(0,R))}^{\vartheta}
 \|V^\infty_{[0,R]}\{e^{it\Delta}f\}\|_{L^{d_n}(B^n(0,R))}^{1-\vartheta}
 \\
 &\lesssim &
 R^{\vartheta b_n/2+\varepsilon}\|f\|_2.
\end{eqnarray*}
The condition $1/r\leq\vartheta/q_n$ gives
$$
 \frac{1/p-1/d_n+(Q_n/r)(1/d_n-1/G_n)}
 {1/2-1/d_n+(Q_n/q_n)(1/d_n-1/G_n)}
 \leq\vartheta.
$$
Thus $\sigma_n(p,r)=A_n(p)+\vartheta b_n$ in this case.

\medskip\noindent
\textbf{Case 2: $d_n\leq p\leq G_n$ and
$\eta_n(p)\leq r\leq\infty$.}
The definition of $\eta_n(p)$ gives
$$
 \frac1p-\frac1{d_n}\leq0,
 \qquad
 \frac1p-\frac1{d_n}
 +\frac{Q_n}{r}\left(\frac1{d_n}-\frac1{G_n}\right)\leq0.
$$
Hence $\sigma_n(p,r)=A_n(p)$, and Lemma \ref{le12}(i) gives
$$
 \|V^r_{[0,R]}\{e^{it\Delta}f\}\|_{L^p(B^n(0,R))}
 \lesssim_\varepsilon R^\varepsilon\|f\|_2.
$$

\medskip\noindent
\textbf{Case 3:}
\[
\left\{
\begin{array}{ll}
	2<p<d_n,
	&
	\displaystyle
	\rho_n(p)\leq r<
	\frac{q_n(1/2-1/d_n)}{1/p-1/d_n};
	\\[10pt]
	d_n\leq p<G_n,
	&
	\rho_n(p)\leq r<\eta_n(p).
\end{array}
\right.
\]
Define
$$
 \lambda=
 \frac{1/p-1/d_n+(Q_n/r)(1/d_n-1/G_n)}
 {1/2-1/d_n+(Q_n/q_n)(1/d_n-1/G_n)},
 \qquad
 \mu=Q_n\left(\frac1r-\frac{\lambda}{q_n}\right),
 \qquad
 \nu=1-\lambda-\mu.
$$
The exclusion of Cases 1 and 2 implies $\lambda>0$ and
$$
 \frac1r\geq
 \frac{1/p-1/d_n}{q_n(1/2-1/d_n)}.
$$
Since $r\geq\rho_n(p)$, we have
$$
\left\{
\begin{aligned}
	\mu
	&=
	\frac{
		Q_n\bigl((1/2-1/d_n)/r-(1/p-1/d_n)/q_n\bigr)
	}{
		1/2-1/d_n+(Q_n/q_n)(1/d_n-1/G_n)
	}
	\ge 0,\\[8pt]
	\nu
	&=
	\frac{
		Q_n(1/2-1/G_n)\bigl(1/\rho_n(p)-1/r\bigr)
	}{
		1/2-1/d_n+(Q_n/q_n)(1/d_n-1/G_n)
	}
	\ge 0.
\end{aligned}
\right.
$$
Moreover, we have 
$$
 \lambda+\mu+\nu=1, 
 \frac1p=\frac\lambda2+\frac\mu{G_n}+\frac\nu{d_n},
 \frac1r=\frac\lambda{q_n}+\frac\mu{Q_n}.
$$
Applying Lemma \ref{le4} twice and then \eqref{eq36},
\begin{eqnarray*}
 \|V^r_{[0,R]}\{e^{it\Delta}f\}\|_{L^p(B^n(0,R))}
 &\leq&
 \|V^{q_n}_{[0,R]}\{e^{it\Delta}f\}\|_{L^2(B^n(0,R))}^{\lambda}
 \|V^{Q_n}_{[0,R]}\{e^{it\Delta}f\}\|_{L^{G_n}(B^n(0,R))}^{\mu}
 \\
 &&\times
 \|V^\infty_{[0,R]}\{e^{it\Delta}f\}\|_{L^{d_n}(B^n(0,R))}^{\nu}
 \\
 &\lesssim &
 R^{\lambda b_n/2+\varepsilon}\|f\|_2.
\end{eqnarray*}
In this case,
$$
 \lambda\geq\max\left\{0,
 \frac{1/p-1/d_n}{1/2-1/d_n}\right\},
 \qquad
 \sigma_n(p,r)=A_n(p)+\lambda b_n,
$$
which implies 
\begin{equation}\label{eq37}
 \|V^r_{[0,R]}\{e^{it\Delta}f\}\|_{L^p(B^n(0,R))}
 \lesssim_\varepsilon
 R^{(\sigma_n(p,r)-A_n(p))/2+\varepsilon}\|f\|_2.
\end{equation}
Applying Theorem \ref{th4} to \eqref{eq37} proves the estimate
for every $s>\sigma_n(p,r)$.
\end{proof}

\section{Proof of the necessity part of Theorem \ref{th1}}
\label{sec5}

\subsection{A lower bound for the variation norm}

\begin{theorem}\label{th6}
Let $n\geq3$, $2\leq p\leq\infty$, $1\leq r\leq\infty$,
and $R\geq1$. Let $f_R\in\mathcal S(\R^n)$ satisfy
$$
 \operatorname{supp}\mathcal F_x(f_R)
 \subset\{\xi\in\R^n:R/2\leq|\xi|\leq2R\}.
$$
Then
\begin{equation}\label{eq38}
 R^{2/r}\|e^{it\Delta}f_R\|_{L_x^pL_t^r([0,1])}
 \lesssim\|V^r_{[0,1]}\{e^{it\Delta}f_R\}\|_{L^p(\R^n)}.
\end{equation}
If \eqref{eq1} holds, then
$$
 \|V^r_{[0,1]}\{e^{it\Delta}f_R\}\|_{L^p(\R^n)}
 \lesssim R^s\|f_R\|_2.
$$
\end{theorem}

\begin{proof}
For $1\leq r<\infty$, choose $0<c<1/64$ and put
$$
 N=\lceil R^2/c\rceil,\qquad h=N^{-1},\qquad
 \frac{c}{1+c}\leq hR^2\leq c.
$$
For $0\leq\tau<h$, the time sequence
$\tau,\tau+h,\ldots,\tau+(N-1)h$ lies in $[0,1]$. Thus
\begin{eqnarray*}
  \int_0^{1-h}
 |e^{i(t+h)\Delta}f_R(x)-e^{it\Delta}f_R(x)|^r\,dt 
 &=&\int_0^h\sum_{j=0}^{N-2}
 |e^{i(\tau+(j+1)h)\Delta}f_R(x)
       -e^{i(\tau+jh)\Delta}f_R(x)|^r\,d\tau\\
 &\leq&h\bigl(V^r_{[0,1]}\{e^{it\Delta}f_R(x)\}\bigr)^r.
\end{eqnarray*}
Taking the spatial norm and changing variables for the backward
difference gives
\begin{eqnarray*}
 && 
 \|e^{i(t+h)\Delta}f_R-e^{it\Delta}f_R\|
       _{L_x^pL_t^r([0,1-h])}}+\|e^{i(t-h)\Delta}f_R-e^{it\Delta}f_R\|
       _{L_x^pL_t^r([h,1])\\
 &&\lesssim h^{1/r}
 \|V^r_{[0,1]}\{e^{it\Delta}f_R\}\|_{L^p(\R^n)}.
\end{eqnarray*}

Choose $\chi\in C_c^\infty(\{\zeta:1/4<|\zeta|<4\})$
with $\chi\equiv1$ on $\{\zeta:1/2\leq|\zeta|\leq2\}$.
On $\operatorname{supp}\chi$,
$$
 0<\frac{c}{16(1+c)}\leq hR^2|\zeta|^2\leq16c<\pi,
 \qquad
 |e^{\pm ihR^2|\zeta|^2}-1|
 \gtrsim hR^2|\zeta|^2\gtrsim c.
$$
For an integer $M>n/2$, it follows that
$$
 \sum_{|\alpha|\leq2M}
 \left\|\partial_\zeta^\alpha
 \left(\frac{\chi(\zeta)}{e^{\pm ihR^2|\zeta|^2}-1}\right)
 \right\|_{L^1_\zeta}
 \lesssim_{M,c}1.
$$
By variable change $\xi=R\zeta$ and integration by parts, we have 
\begin{eqnarray*}
 \left|\mathcal F_x^{-1}
 \left(\frac{\chi(\xi/R)}{e^{\pm ih|\xi|^2}-1}\right)(y)\right|
 &\lesssim_{M,c}&R^n(1+R^2|y|^2)^{-M}\end{eqnarray*}
which indicates 
 \begin{eqnarray*}
 \left\|\mathcal F_x^{-1}
 \left(\frac{\chi(\xi/R)}{e^{\pm ih|\xi|^2}-1}\right)\right\|_{L^1}
 &\lesssim_{M,c}&R^n\int_{\R^n}(1+R^2|y|^2)^{-M}\,dy
 \lesssim_{M,c}1.
\end{eqnarray*}
Since
$$
 \mathcal F_x\bigl(e^{i(t\pm h)\Delta}f_R-e^{it\Delta}f_R\bigr)(\xi)
 =\bigl(e^{\mp ih|\xi|^2}-1\bigr)
 e^{-it|\xi|^2}\mathcal F_x(f_R)(\xi),
$$
Minkowski's inequality and Young's inequality yield
\begin{eqnarray*}
\|e^{it\Delta}f_R\|_{L_x^pL_t^r([0,1-h])}
 &\leq&
 \left\|\mathcal F_x^{-1}
 \left(\frac{\chi(\xi/R)}{e^{-ih|\xi|^2}-1}\right)\right\|_{L^1}\|e^{i(t+h)\Delta}f_R-e^{it\Delta}f_R\|
 _{L_x^pL_t^r([0,1-h])}\\
 &\lesssim&h^{1/r}
 \|V^r_{[0,1]}\{e^{it\Delta}f_R\}\|_{L^p}
\end{eqnarray*}
and
\begin{eqnarray*}
 \|e^{it\Delta}f_R\|_{L_x^pL_t^r([h,1])}&\leq& \left\|\mathcal F_x^{-1}
 \left(\frac{\chi(\xi/R)}{e^{ih|\xi|^2}-1}\right)\right\|_{L^1}  \|e^{i(t-h)\Delta}f_R-e^{it\Delta}f_R\|
 _{L_x^pL_t^r([h,1])}\\
 &\lesssim&h^{1/r}
 \|V^r_{[0,1]}\{e^{it\Delta}f_R\}\|_{L^p}.
\end{eqnarray*}
As $h<1/2$ and $[0,1-h]\cup[h,1]=[0,1]$,
\begin{eqnarray*}
 \|e^{it\Delta}f_R\|_{L_x^pL_t^r([0,1])}
 &\leq&\|e^{it\Delta}f_R\|_{L_x^pL_t^r([0,1-h])}
       +\|e^{it\Delta}f_R\|_{L_x^pL_t^r([h,1])}\\
 &\lesssim&R^{-2/r}
 \|V^r_{[0,1]}\{e^{it\Delta}f_R\}\|_{L^p}.
\end{eqnarray*}
For $r=\infty$, \eqref{eq38} follows from the definition of $V^\infty$.
If \eqref{eq1} holds, then
\begin{eqnarray*}
 \|V^r_{[0,1]}\{e^{it\Delta}f_R\}\|_{L^p(\R^n)}
 \lesssim\|f_R\|_{H^s(\R^n)}\lesssim R^s\|f_R\|_2.
\end{eqnarray*}
\end{proof}

\subsection{Counterexample in the range \texorpdfstring{$s<A_n(p)$}{s < A n(p)}}

\begin{theorem}\label{th7}
Let $n\geq3$, $2\leq p\leq\infty$, and $1\leq r\leq\infty$.
If
$$
 \bigl\|V^r_{[0,1]}\{e^{it\Delta}f\}\bigr\|_{L^p(\R^n)}
 \lesssim\|f\|_{H^s(\R^n)}
$$
holds for all $f\in H^s(\R^n)$, then necessarily
$$
 s\geq A_n(p).
$$
\end{theorem}

\begin{proof}
Choose $0\leq a\in C_c^\infty(B^n(e_1,1/8))$, $a\not\equiv0$,
and define
$$
 \mathcal F_x(f_R)(\xi)=R^{-n/2}a(\xi/R).
$$
By variable change $\xi=R\zeta$,
$$
 \|f_R\|_2^2\approx\int_{\R^n}|a(\zeta)|^2\,d\zeta\approx1,
 \qquad \|f_R\|_{H^s}\lesssim R^s.
$$
For $|x|\leq cR^{-1}$, $0\leq t\leq cR^{-2}$, and
$\xi\in\operatorname{supp}\mathcal F_x(f_R)$,
$$
 |x\cdot\xi-t|\xi|^2|
 \leq |x||\xi|+t|\xi|^2
 \leq2c+4c\leq\frac\pi3,
$$
provided that $c>0$ is sufficiently small. Hence
\begin{eqnarray*}
 |e^{it\Delta}f_R(x)|
 \geq\operatorname{Re}\bigl(e^{it\Delta}f_R(x)\bigr)
 \geq\frac{R^{-n/2}}{2(2\pi)^n}
 \int_{\R^n}a(\xi/R)\,d\xi
 \gtrsim R^{n/2}.
\end{eqnarray*}
For $1\leq r<\infty$, Theorem \ref{th6} gives
\begin{eqnarray*}
 \|V^r_{[0,1]}\{e^{it\Delta}f_R\}\|_{L^p}
 \gtrsim R^{2/r}\|e^{it\Delta}f_R\|_{L_x^pL_t^r([0,1])}
 \gtrsim R^{2/r}R^{n/2}(cR^{-2})^{1/r}
               |B^n(0,cR^{-1})|^{1/p}\gtrsim R^{n/2-n/p}.
\end{eqnarray*}
For $r=\infty$, taking $t=0$ gives
$$
 \|V^\infty_{[0,1]}\{e^{it\Delta}f_R\}\|_{L^p}
 \geq2\|f_R\|_{L^p(B^n(0,cR^{-1}))}
 \gtrsim R^{n/2-n/p}.
$$
Thus
$$
 R^{A_n(p)}
 \lesssim\|V^r_{[0,1]}\{e^{it\Delta}f_R\}\|_{L^p}
 \lesssim\|f_R\|_{H^s}
 \lesssim R^s,
$$
which requires $s\geq A_n(p)$ as $R\to\infty$.
\end{proof}

\subsection{Counterexample in the range \texorpdfstring{$s<1/p+1/r$}{s < 1/p+1/r}}

\begin{theorem}\label{th8}
Let $n\geq3$, $2\leq p\leq\infty$, and $1\leq r\leq\infty$.
If
$$
 \bigl\|V^r_{[0,1]}\{e^{it\Delta}f\}\bigr\|_{L^p(\R^n)}
 \lesssim\|f\|_{H^s(\R^n)}
$$
holds for all $f\in H^s(\R^n)$, then necessarily
$$
 s\geq \frac1p+\frac1r.
$$
\end{theorem}

\begin{proof}
Choose $0\leq\widehat\phi\in C_c^\infty(B^n(0,\delta))$,
$\widehat\phi\not\equiv0$, with $0<\delta<1/100$, and set
$$
 f_R(x)=e^{iRx_1}\phi(x).
$$
Then
$$
 \operatorname{supp}\mathcal F_x(f_R)\subset B^n(Re_1,\delta),
 \qquad \|f_R\|_2=\|\phi\|_2,\qquad
 \|f_R\|_{H^s}\lesssim R^s\|\phi\|_2.
$$
By variable change $\xi=Re_1+\eta$,
$$
 e^{it\Delta}f_R(x)
 =e^{i(Rx_1-R^2t)}(e^{it\Delta}\phi)(x-2Rte_1).
$$
Let
$$
 \Omega_R=\{x\in\R^n:R/2\leq x_1\leq R,\ |x'|\leq1/4\},
 \qquad |\Omega_R|\approx R.
$$
For $x\in\Omega_R$, put
$$
 t_x=\frac{x_1}{2R},\qquad
 I_x=\left[t_x-\frac1{8R},t_x+\frac1{8R}\right]\subset(0,1).
$$
If $t\in I_x$, then
$$
 |x_1-2Rt|\leq\frac14,\qquad
 |x-2Rte_1|\leq\frac12,
$$
and, for $\eta\in\operatorname{supp}\widehat\phi$,
$$
 |(x-2Rte_1)\cdot\eta-t|\eta|^2|
 \leq\frac\delta2+\delta^2\leq\frac\pi3.
$$
Consequently,
$$
 \operatorname{Re}\bigl((e^{it\Delta}\phi)(x-2Rte_1)\bigr)
 \geq\frac1{2(2\pi)^n}\int_{\R^n}\widehat\phi(\eta)\,d\eta
 \gtrsim1.
$$
For finite $r$, choose
$$
 t_j=t_x-\frac1{8R}+\frac{j\pi}{R^2},\qquad
 0\leq j\leq J_R=\left\lfloor\frac R{4\pi}\right\rfloor.
$$
For sufficiently large $R$,
$$
 t_x-\frac1{8R}\leq t_0<\cdots<t_{J_R}
 \leq t_x+\frac1{8R},\qquad J_R\gtrsim R,
 \qquad e^{-iR^2t_{j+1}}=-e^{-iR^2t_j}.
$$
It follows that, for $0\leq j<J_R$,
\begin{eqnarray*}
 |e^{it_{j+1}\Delta}f_R(x)-e^{it_j\Delta}f_R(x)|
 &\geq&\operatorname{Re}\bigl((e^{it_{j+1}\Delta}\phi)
                 (x-2Rt_{j+1}e_1)\bigr)+\operatorname{Re}\bigl((e^{it_j\Delta}\phi)
                 (x-2Rt_je_1)\bigr) \gtrsim1.
\end{eqnarray*}
Thus, for $x\in\Omega_R$,
\begin{eqnarray*}
 V^r_{[0,1]}\{e^{it\Delta}f_R(x)\}
 \geq\left(\sum_{j=0}^{J_R-1}
 |e^{it_{j+1}\Delta}f_R(x)-e^{it_j\Delta}f_R(x)|^r\right)^{1/r}
 \gtrsim J_R^{1/r}\gtrsim R^{1/r}.
\end{eqnarray*}
For $r=\infty$, taking $t=t_x$ gives
$$
 V^\infty_{[0,1]}\{e^{it\Delta}f_R(x)\}
 \geq2|(e^{it_x\Delta}\phi)(x-2Rt_xe_1)|\gtrsim1.
$$
Therefore, with $1/\infty=0$,
\begin{eqnarray*}
 R^{1/p+1/r}
 &\lesssim&R^{1/r}|\Omega_R|^{1/p}
 \lesssim\|V^r_{[0,1]}\{e^{it\Delta}f_R\}\|_{L^p(\Omega_R)}
 \lesssim\|f_R\|_{H^s}\lesssim R^s.
\end{eqnarray*}
Letting $R\to\infty$ yields $s\geq1/p+1/r$.
\end{proof}

\subsection{Counterexample in the range \texorpdfstring{$s<2/r+1/p-1/2$}{s < 2/r+1/p-1/2}}

\begin{theorem}\label{th9}
Let $n\geq3$, $2\leq p\leq\infty$, and $1\leq r<\infty$.
If
$$
 \bigl\|V^r_{[0,1]}\{e^{it\Delta}f\}\bigr\|_{L^p(\R^n)}
 \lesssim\|f\|_{H^s(\R^n)}
$$
holds for all $f\in H^s(\R^n)$, then necessarily
$$
 s\geq \frac2r+\frac1p-\frac12.
$$
\end{theorem}

\begin{proof}
Choose $0\leq a\in C_c^\infty(B^n(0,\delta))$,
$a\not\equiv0$, with $0<\delta<1/100$, and define
$$
 \mathcal F_x(f_R)(\xi)
 =R^{1/2}a\bigl(R(\xi_1-R),\xi'\bigr).
$$
By variable change $\eta_1=R(\xi_1-R)$, $\eta'=\xi'$,
$$
 \|f_R\|_2^2\approx\int_{\R^n}|a(\eta)|^2\,d\eta\approx1,
 \qquad \|f_R\|_{H^s}\lesssim R^s,
$$
and
\begin{eqnarray*}
 e^{it\Delta}f_R(x)
 &=&\frac{R^{-1/2}e^{i(Rx_1-R^2t)}}{(2\pi)^n}\int_{\R^n}
 e^{i\left((x_1/R-2t)\eta_1+x'\cdot\eta'
             -t(R^{-2}\eta_1^2+|\eta'|^2)\right)}a(\eta)\,d\eta.
\end{eqnarray*}
Let
$$
 \Omega_R=\{x\in\R^n:|x_1|\leq R/4,\ |x'|\leq1/4\},
 \qquad |\Omega_R|\approx R.
$$
For $x\in\Omega_R$, $0\leq t\leq1$, and
$\eta\in\operatorname{supp}a$,
\begin{eqnarray*}
 \lefteqn{\left|(x_1/R-2t)\eta_1+x'\cdot\eta'
                   -t(R^{-2}\eta_1^2+|\eta'|^2)\right|}\\
 &\leq&\left(\frac14+2\right)\delta
                  +\frac\delta4+2\delta^2
 \leq\frac\pi3.
\end{eqnarray*}
Hence
$$
 \operatorname{Re}\bigl(e^{-i(Rx_1-R^2t)}e^{it\Delta}f_R(x)\bigr)
 \geq\frac{R^{-1/2}}{2(2\pi)^n}\int_{\R^n}a(\eta)\,d\eta
 \gtrsim R^{-1/2}.
$$
Take
$$
 t_j=\frac{j\pi}{R^2},\qquad
 0\leq j\leq J_R=\left\lfloor\frac{R^2}{\pi}\right\rfloor.
$$
Then $0\leq t_0<\cdots<t_{J_R}\leq1$, $J_R\gtrsim R^2$,
and $e^{-iR^2t_{j+1}}=-e^{-iR^2t_j}$. Therefore,
\begin{eqnarray*}
|e^{it_{j+1}\Delta}f_R(x)-e^{it_j\Delta}f_R(x)|
 \geq\operatorname{Re}\bigl(e^{-i(Rx_1-R^2t_{j+1})} e^{it_{j+1}\Delta}f_R(x)\bigr)+\operatorname{Re}\bigl(e^{-i(Rx_1-R^2t_j)}e^{it_j\Delta}f_R(x)\bigr)\gtrsim R^{-1/2}.
\end{eqnarray*}
For  $x\in\Omega_R$, we have 
\begin{eqnarray*}
 V^r_{[0,1]}\{e^{it\Delta}f_R(x)\}
 \geq\left(\sum_{j=0}^{J_R-1}
 |e^{it_{j+1}\Delta}f_R(x)-e^{it_j\Delta}f_R(x)|^r\right)^{1/r}
 \gtrsim R^{-1/2}J_R^{1/r}
 \gtrsim R^{2/r-1/2},
\end{eqnarray*}
which implies 
\begin{eqnarray*}
 R^{2/r+1/p-1/2}\lesssim R^{2/r-1/2}|\Omega_R|^{1/p}\lesssim\|V^r_{[0,1]}\{e^{it\Delta}f_R\}\|_{L^p(\Omega_R)}
 \lesssim \|f_R\|_{H^s}\lesssim R^s.
\end{eqnarray*}
Letting $R\to\infty$ yields $s\geq2/r+1/p-1/2$.
\end{proof}

\subsection{Counterexample in the range \texorpdfstring{$s<A_n(p)+2/r-1$}{s < A n(p)+2/r-1}}

\begin{theorem}\label{th10}
Let $n\geq3$, $2\leq p\leq\infty$, and $1\leq r\leq\infty$.
If
$$
 \bigl\|V^r_{[0,1]}\{e^{it\Delta}f\}\bigr\|_{L^p(\R^n)}
 \lesssim\|f\|_{H^s(\R^n)}
$$
holds for all $f\in H^s(\R^n)$, then necessarily
$$
 s\geq A_n(p)+\frac2r-1.
$$
\end{theorem}

\begin{proof}
Choose $0\leq a\in C_c^\infty((0,\delta))$, $a\not\equiv0$,
with $0<\delta<1/100$, and set
$$
 \mathcal F_x(f_R)(\xi)
 =R^{-(n-2)/2}a\bigl(R(|\xi|-R)\bigr).
$$
Then
$$
 \operatorname{supp}\mathcal F_x(f_R)
 \subset\{\xi\in\R^n:R\leq|\xi|\leq R+\delta/R\}.
$$
By polar coordinates and variable change $v=R(\rho-R)$,
$$
 \|f_R\|_2^2
 \approx\int_0^\delta|a(v)|^2
            \left(1+\frac v{R^2}\right)^{n-1}\,dv
 \approx1,\qquad \|f_R\|_{H^s}\lesssim R^s.
$$
On the Fourier support,
$$
 0\leq|\xi|^2-R^2
 \leq\frac\delta R\left(2R+\frac\delta R\right)
 \leq3\delta.
$$
For $|x|\leq cR^{-1}$ and $0\leq t\leq1$,
$$
 |x\cdot\xi-t(|\xi|^2-R^2)|
 \leq |x||\xi|+t\bigl||\xi|^2-R^2\bigr|
 \leq2c+3\delta\leq\frac\pi3
$$
if $c>0$ is sufficiently small. Therefore,
\begin{eqnarray*}
 \operatorname{Re}\bigl(e^{itR^2}e^{it\Delta}f_R(x)\bigr)
 &\geq&\frac1{2(2\pi)^n}\int_{\R^n}\mathcal F_x(f_R)(\xi)\,d\xi\\
 &\gtrsim&R^{-(n-2)/2}
   \int_R^{R+\delta/R}a\bigl(R(\rho-R)\bigr)\rho^{n-1}\,d\rho\\
 &\gtrsim&R^{(n-2)/2}\int_0^\delta a(v)\,dv
 \gtrsim R^{(n-2)/2}.
\end{eqnarray*}
For finite $r$, take
$$
 t_j=\frac{j\pi}{R^2},\qquad
 0\leq j\leq J_R=\left\lfloor\frac{R^2}{\pi}\right\rfloor.
$$
Since $J_R\gtrsim R^2$ and $e^{it_jR^2}=(-1)^j$,
\begin{eqnarray*}
  |e^{it_{j+1}\Delta}f_R(x)-e^{it_j\Delta}f_R(x)| 
 \geq\operatorname{Re}\bigl(e^{it_{j+1}R^2}e^{it_{j+1}\Delta}f_R(x)\bigr) +\operatorname{Re}\bigl(e^{it_jR^2}
                              e^{it_j\Delta}f_R(x)\bigr)
 \gtrsim R^{(n-2)/2}.
\end{eqnarray*}
Consequently, for $|x|\leq cR^{-1}$,
\begin{eqnarray*}
 V^r_{[0,1]}\{e^{it\Delta}f_R(x)\}
 &\geq&\left(\sum_{j=0}^{J_R-1}
 |e^{it_{j+1}\Delta}f_R(x)-e^{it_j\Delta}f_R(x)|^r\right)^{1/r}\\
 &\gtrsim&R^{(n-2)/2}J_R^{1/r}\\
 &\gtrsim& R^{(n-2)/2+2/r}.
\end{eqnarray*}
Thus
\begin{eqnarray*}
 R^{A_n(p)+2/r-1}
 &\lesssim&R^{(n-2)/2+2/r}|B^n(0,cR^{-1})|^{1/p}\\
 &\lesssim&\|V^r_{[0,1]}\{e^{it\Delta}f_R\}\|_{L^p}
 \lesssim\|f_R\|_{H^s}\lesssim R^s.
\end{eqnarray*}
For $r=\infty$, taking $t=0$ gives
$$
 \|V^\infty_{[0,1]}\{e^{it\Delta}f_R\}\|_{L^p}
 \geq2\|f_R\|_{L^p(B^n(0,cR^{-1}))}
 \gtrsim R^{(n-2)/2-n/p}.
$$
Letting $R\to\infty$ proves $s\geq A_n(p)+2/r-1$.
\end{proof}

\subsection{The lower bound at \texorpdfstring{$p=2$}{p=2}}

\begin{theorem}\label{th11}
Let $n\geq3$ and $1\leq r\leq\infty$.
If
$$
 \bigl\|V^r_{[0,1]}\{e^{it\Delta}f\}\bigr\|_{L^2(\R^n)}
 \lesssim\|f\|_{H^s(\R^n)}
$$
holds for all $f\in H^s(\R^n)$, then necessarily
$$
 s\geq b_n.
$$
\end{theorem}

\begin{proof}
Fix $0<\delta<b_n/2$. Bourgain's construction \cite{B2016},
as detailed in \cite[(2.10)--(2.11) and Section 7]{Pierce2020},
gives, for arbitrarily large $L$, nonzero $h_L\in\mathcal S(\R^n)$
such that
$$
 \operatorname{supp}\mathcal F_x(h_L)
 \subset\{\xi\in\R^n:C^{-1}L\leq|\xi|\leq CL\},
$$
and
$$
 \left\|\sup_{0<t<c/L}|e^{it\Delta}h_L|\right\|_{L^1(B^n(0,1))}
 \gtrsim_\delta L^{n/(2(n+1))-\delta}\|h_L\|_2,
$$
where $C\geq1$ and $0<c<1$ are fixed. By H\"older's inequality,
\begin{eqnarray}\label{eq39}
 \left\|\sup_{0<t<c/L}|e^{it\Delta}h_L|\right\|_{L^2(B^n(0,1))}
 &\geq&|B^n(0,1)|^{-1/2}
 \left\|\sup_{0<t<c/L}|e^{it\Delta}h_L|\right\|_{L^1(B^n(0,1))}
 \nonumber\\
 &\gtrsim_\delta&L^{n/(2(n+1))-\delta}\|h_L\|_2.
\end{eqnarray}
Set $L=R^2$ and
$$
 f_R(x)=R^{-n/2}h_{R^2}(x/R).
$$
By variable change,
$$
 \mathcal F_x(f_R)(\xi)=R^{n/2}\mathcal F_x(h_{R^2})(R\xi),
 \qquad \|f_R\|_2=\|h_{R^2}\|_2,
$$
so
$$
 \operatorname{supp}\mathcal F_x(f_R)
 \subset\{\xi\in\R^n:C^{-1}R\leq|\xi|\leq CR\},
 \qquad \|f_R\|_{H^s}\lesssim R^s\|f_R\|_2.
$$
Since
$$
 e^{it\Delta}f_R(x)
 =R^{-n/2}(e^{i(t/R^2)\Delta}h_{R^2})(x/R),
$$
\eqref{eq39} and variable change $x=Ry$, $t=R^2\tau$ give
\begin{eqnarray*}
 \left\|\sup_{0<t<c}|e^{it\Delta}f_R|\right\|_{L^2(\R^n)}
 &\geq&\left\|\sup_{0<t<c}|e^{it\Delta}f_R|\right\|_{L^2(B^n(0,R))}\\
 &=&\left\|\sup_{0<\tau<c/R^2}|e^{i\tau\Delta}h_{R^2}|\right\|
       _{L^2(B^n(0,1))}\\
 &\gtrsim_\delta&R^{b_n-2\delta}\|f_R\|_2.
\end{eqnarray*}
For finite $r$,
$$
 \sup_{0<t<c}|e^{it\Delta}f_R(x)|
 \leq |f_R(x)|+V^r_{[0,1]}\{e^{it\Delta}f_R(x)\}.
$$
Hence, for all sufficiently large $R$,
\begin{eqnarray*}
 \|V^r_{[0,1]}\{e^{it\Delta}f_R\}\|_{L^2}
 &\geq&\left\|\sup_{0<t<c}|e^{it\Delta}f_R|\right\|_{L^2}
       -\|f_R\|_2\\
 &\geq&(c_\delta R^{b_n-2\delta}-1)\|f_R\|_2\\
 &\geq&\frac{c_\delta}{2}R^{b_n-2\delta}\|f_R\|_2.
\end{eqnarray*}
For $r=\infty$,
$$
 \|V^\infty_{[0,1]}\{e^{it\Delta}f_R\}\|_{L^2}
 \geq2\left\|\sup_{0<t<c}|e^{it\Delta}f_R|\right\|_{L^2}
 \gtrsim_\delta R^{b_n-2\delta}\|f_R\|_2.
$$
In either case,
$$
 R^{b_n-2\delta}\|f_R\|_2
 \lesssim_\delta\|V^r_{[0,1]}\{e^{it\Delta}f_R\}\|_{L^2}
 \lesssim\|f_R\|_{H^s}
 \lesssim R^s\|f_R\|_2.
$$
Since $\|f_R\|_2>0$, letting $R\to\infty$ gives
$s\geq b_n-2\delta$. Letting $\delta\downarrow0$ proves
$s\geq b_n$.
\end{proof}

The first four examples give all the lower bounds in
parts (i)--(iv) of Theorem \ref{th1}; the last example gives
part (v). Together with Section \ref{sec4}, this
completes the proof of the theorem.

\noindent\textbf{Data availability} This manuscript no associate data.

\noindent\textbf{A conflict of interest statement} This manuscript is original, has not been formally published,and is not under consideration for publication elsewhere. Both authors have approved this submission. The authors declare no competing interests.

\noindent\textbf{Statement on the use of AI.}
The authors developed the main ideas, wrote the proofs, and
determined the overall direction of this work.
We suggested to GPT-6 Astra that the two-dimensional counterexample in Chen and Tang \cite{CT2026} could be used to obtain a counterexample for higher-dimensional Schr\"odinger variation estimates, and asked the model to help work out the details of this extension in Section 5. We also used the model to refine the wording of the manuscript. Figures 1 is generated by GPT-6 Astra based on our descriptions of their intended meaning and scope. The proofs for Section 3  in this paper were completed by the authors. We take full responsibility for the mathematical content of
the paper.


\begin{thebibliography}{99}
	
	\bibitem{B1989}
	J. Bourgain,
	\emph{Pointwise ergodic theorems for arithmetic sets},
	Publ. Math. Inst. Hautes \`Etudes Sci. \textbf{69} (1989), 5--45.
	
	\bibitem{B1995}
	J. Bourgain,
	\emph{Some new estimates on oscillatory integrals},
	in \emph{Essays on Fourier Analysis in Honor of Elias M. Stein},
	Princeton Math. Ser., vol. 42,
	Princeton University Press, Princeton, NJ, 1995, pp. 83--112.
	
	\bibitem{B2013}
	J. Bourgain,
	\emph{On the Schr\"odinger maximal function in higher dimension},
	Proc. Steklov Inst. Math. \textbf{280} (2013), no. 1, 46--60.
	
	\bibitem{B2016}
	J. Bourgain,
	\emph{A note on the Schr\"odinger maximal function},
	J. Anal. Math. \textbf{130} (2016), no. 1, 393--396.
	
	\bibitem{BD2015}
	J. Bourgain and C. Demeter,
	\emph{The proof of the $\ell^2$ decoupling conjecture},
	Ann. of Math. (2) \textbf{182} (2015), no. 1, 351--389.
	
	\bibitem{CJRW2000}
	J. T. Campbell, R. L. Jones, K. Reinhold, and M. Wierdl,
	\emph{Oscillation and variation for the Hilbert transform},
	Duke Math. J. \textbf{105} (2000), no. 1, 59--83.
	
	\bibitem{CMW2021}
	Z. Cao, C. Miao, and M. Wang,
	\emph{$L^p$-estimate of Schr\"odinger maximal function in higher dimensions},
	J. Funct. Anal. \textbf{281} (2021), no. 5, Art. 109091.
	
	\bibitem{C1980}
	L. Carleson,
	\emph{Some analytic problems related to statistical mechanics},
	in \emph{Euclidean Harmonic Analysis (Proc. Sem., Univ. Maryland,
		College Park, Md., 1979)}, Lecture Notes in Math., vol. 779,
	Springer, Berlin, 1980, pp. 5--45.
	
	\bibitem{CT2026}
	Y. Chen and Y. Tang,
	\emph{Sharp variational estimates for Schr\"odinger equation in $\R^2$},
	arXiv:2610.31099 [math.CA], 2026.
	\newblock \url{https://arxiv.org/abs/2610.31099}.
	
	\bibitem{DK1982}
	B. E. J. Dahlberg and C. E. Kenig,
	\emph{A note on the almost everywhere behavior of solutions to the
		Schr\"odinger equation},
	in \emph{Harmonic Analysis (Minneapolis, Minn., 1981)},
	Lecture Notes in Math., vol. 908, Springer, Berlin--New York, 1982,
	pp. 205--209.
	
	\bibitem{DG2016}
	C. Demeter and S. Guo,
	\emph{Schr\"odinger maximal function estimates via the pseudoconformal transformation},
	arXiv:1608.07640 [math.CA], 2016.
	
	\bibitem{DGL2017}
	X. Du, L. Guth, and X. Li,
	\emph{A sharp Schr\"odinger maximal estimate in $\mathbb{R}^2$},
	Ann. of Math. (2) \textbf{186} (2017), no. 2, 607--640.
	
	\bibitem{DKWZ2020}
	X. Du, J. Kim, H. Wang, and R. Zhang,
	\emph{Lower bounds for estimates of the Schr\"odinger maximal function},
	Math. Res. Lett. \textbf{27} (2020), no. 3, 687--692.
	
	\bibitem{DL2025}
	X. Du and J. Li,
	\emph{$L^p$ estimates of the maximal Schr\"odinger operator in $\R^n$},
	J. Funct. Anal. \textbf{288} (2025), no. 3, Art. 110737.
	
	\bibitem{DZ2019}
	X. Du and R. Zhang,
	\emph{Sharp $L^2$ estimates of the Schr\"odinger maximal function in
		higher dimensions},
	Ann. of Math. (2) \textbf{189} (2019), no. 3, 837--861.
	
	\bibitem{GRY2020}
	S. Guo, J. Roos, and P.-L. Yung,
	\emph{Sharp variation-norm estimates for oscillatory integrals related
		to Carleson's theorem},
	Anal. PDE \textbf{13} (2020), no. 5, 1457--1500.
	
	\bibitem{JSW2008}
	R. L. Jones, A. Seeger, and J. Wright,
	\emph{Strong variational and jump inequalities in harmonic analysis},
	Trans. Amer. Math. Soc. \textbf{360} (2008), no. 12, 6711--6742.
	
	\bibitem{KT1998}
	M. Keel and T. Tao,
	\emph{Endpoint Strichartz estimates},
	Amer. J. Math. \textbf{120} (1998), no. 5, 955--980.
	
	\bibitem{L2006}
	S. Lee,
	\emph{On pointwise convergence of the solutions to Schr\"odinger
		equations in $\mathbb{R}^2$},
	Int. Math. Res. Not. \textbf{2006} (2006), Art. ID 32597, 21 pp.
	
	\bibitem{LRV2011}
	S. Lee, K. M. Rogers, and A. Vargas,
	\emph{An endpoint space--time estimate for the Schr\"odinger equation},
	Adv. Math. \textbf{226} (2011), no. 5, 4266--4285.
	
	\bibitem{LRS2013}
	S. Lee, K. M. Rogers, and A. Seeger,
	\emph{On space--time estimates for the Schr\"odinger operator},
	J. Math. Pures Appl. (9) \textbf{99} (2013), no. 1, 62--85.
	
	\bibitem{L1976}
	D. L\'epingle,
	\emph{La variation d'ordre $p$ des semi-martingales},
	Z. Wahrscheinlichkeitstheorie verw. Gebiete \textbf{36} (1976), no. 4,
	295--316.
	
	\bibitem{LMY2024}
	J. Li, C. Miao, and A. Yu,
	\emph{The space-time estimates for the Schr\"odinger equation},
	arXiv:2402.13539 [math.CA], 2024.
	
	\bibitem{LR2017}
	R. Luc\`a and K. M. Rogers,
	\emph{Coherence on fractals versus pointwise convergence for the
		Schr\"odinger equation},
	Comm. Math. Phys. \textbf{351} (2017), no. 1, 341--359.
	
	\bibitem{MVV1996}
	A. Moyua, A. Vargas, and L. Vega,
	\emph{Schr\"odinger maximal function and restriction properties of the Fourier transform},
	International Mathematics Research Notices \textbf{1996} (1996),
	no. 16, 793--815.
	
	\bibitem{OSTTW2012}
	R. Oberlin, A. Seeger, T. Tao, C. Thiele, and J. Wright,
	\emph{A variation norm Carleson theorem},
	J. Eur. Math. Soc. (JEMS) \textbf{14} (2012), no. 2, 421--464.
	
	\bibitem{Pierce2020}
	L. B. Pierce,
	\emph{On Bourgain's counterexample for the Schr\"odinger maximal function},
	Quart. J. Math. \textbf{71} (2020), no. 4, 1309--1344.
	
	\bibitem{PX1988}
	G. Pisier and Q. Xu,
	\emph{The strong $p$-variation of martingales and orthogonal series},
	Probab. Theory Related Fields \textbf{77} (1988), no. 4, 497--514.
	
	\bibitem{R2008}
	K. M. Rogers,
	\emph{A local smoothing estimate for the Schr\"odinger equation},
	Adv. Math. \textbf{219} (2008), no. 6, 2105--2122.
	
	\bibitem{RS2010}
	K. M. Rogers and A. Seeger,
	\emph{Endpoint maximal and smoothing estimates for Schr\"odinger equations},
	J. Reine Angew. Math. \textbf{640} (2010), 47--66.
	
	\bibitem{RV2007}
	K. M. Rogers and P. Villarroya,
	Global estimates for the Schr\"odinger maximal operator,
	Ann. Acad. Sci. Fenn. Math. \textbf{32} (2007), 425--435.
	
	\bibitem{S1987}
	P. Sj\"olin,
	\emph{Regularity of solutions to the Schr\"odinger equation},
	Duke Math. J. \textbf{55} (1987), no. 3, 699--715.
	
	\bibitem{Stein1993}
	E. M. Stein,
	\emph{Harmonic Analysis: Real-Variable Methods, Orthogonality, and Oscillatory Integrals},
	Princeton Mathematical Series, vol. 43,
	Princeton University Press, Princeton, NJ, 1993.
	
	\bibitem{SW2001}
	E. M. Stein and S. Wainger,
	\emph{Oscillatory integrals related to Carleson's theorem},
	Math. Res. Lett. \textbf{8} (2001), nos. 5--6, 789--800.
	
	\bibitem{Tao2003}
	T. Tao,
	\emph{A sharp bilinear restriction estimate for paraboloids},
	Geom. Funct. Anal. \textbf{13} (2003), no. 6, 1359--1384.
	
	\bibitem{TV2000}
	T. Tao and A. Vargas,
	\emph{A bilinear approach to cone multipliers. II. Applications},
	Geom. Funct. Anal. \textbf{10} (2000), no. 1, 216--258.
	
	\bibitem{Tomas1975}
	P. A. Tomas,
	\emph{A restriction theorem for the Fourier transform},
	Bull. Amer. Math. Soc. \textbf{81} (1975), no. 2, 477--478.
	
	\bibitem{V1988}
	L. Vega,
	\emph{Schr\"odinger equations: pointwise convergence to the initial data},
	Proc. Amer. Math. Soc. \textbf{102} (1988), no. 4, 874--878.
	
\end{thebibliography}
\end{document}